\documentclass[11pt]{amsart}

\usepackage{amsthm,amssymb,amsfonts,amsmath,enumerate,graphicx,color,tikz}
 \usepackage{mathtools}

\usepackage{mathrsfs}
\usepackage[margin=1in]{geometry}

\usepackage{caption}
\usepackage{subcaption}

\usepackage[utf8]{inputenc}
\usepackage{palatino}

\usepackage{hyperref}
	\hypersetup{
	colorlinks=true,
	linkcolor={black},
	citecolor={blue}
	}
	
	\usepackage{tikz}
\usetikzlibrary{calc}
\usepackage{pgfplots}
\pgfplotsset{compat=1.11}
\usepackage{tikz-3dplot}

\newcommand\inner[1]{\langle {#1} \rangle}%
\newcommand\y{\mathbf{y}}%
\newcommand\OV[1]{\overline{#1}}%
\newcommand{\vol}{\operatorname{vol}}%
\newcommand{\Vol}{\operatorname{Vol}}%
\newcommand\Vdown{V^{\downarrow}}%
\newcommand\convDown[1]{\{#1\}^{\downarrow}}%
\newcommand\U{\mathsf{U}}%
\newcommand\PM{\{-1,+1\}}%
\newcommand\Triang{\mathcal{T}}%
\newcommand\Subdiv{\mathcal{S}}%
\newcommand\1{\mathbf{1}}%
\newcommand\w{\mathbf{w}}%
\newcommand{\Mat}{\mathsf{Mat}}%
\newcommand\Chain{\mathcal{C}}%
\newcommand\Anti{\mathcal{A}}%
\newcommand\UC{\mathsf{U}\Chain}%
\newcommand{\ehr}{\operatorname{ehr}}
\renewcommand{\k}{\Bbbk}

\newcommand{\C}{\mathbb{C}}
\newcommand{\Z}{\mathbb{Z}}
\newcommand{\R}{\mathbb{R}}
\newcommand{\Q}{\mathbb{Q}}
\newcommand{\conv}{\mathrm{conv}}
\renewcommand{\emptyset}{\varnothing}
\newcommand{\Zpos}{\Z_{\geq 1}}

\renewcommand\a{{\boldsymbol a}}
\renewcommand\c{{\boldsymbol c}}
\renewcommand\d{{\boldsymbol d}}
\newcommand\e{{\boldsymbol e}}

\newcommand\p{{\boldsymbol p}}
\newcommand\q{{\boldsymbol q}}
\newcommand\s{{\boldsymbol s}}
\renewcommand\r{{\boldsymbol r}}
\renewcommand\t{{\boldsymbol t}}
\renewcommand\u{{\boldsymbol u}}
\newcommand\vb{{\boldsymbol v}}
\newcommand\x{{\boldsymbol x}}
\newcommand\z{{\boldsymbol z}}

\newcommand\ceq{\coloneqq}

\newcommand\Def[1]{\textbf{#1}}%

\newcommand{\rem}[1]{}

\newcommand{\BB}{\mathcal{BB}} 
\newcommand{\pos}{\mathcal{BB}_+} 
\newcommand{\Cn}{\mathcal{C}(n)}
\newcommand{\Cposn}{\mathcal{C}_+(n)}

\newcommand{\sgn}{\operatorname{sgn}}

\newtheorem{theorem}{Theorem}[section]
\newtheorem{corollary}[theorem]{Corollary}
\newtheorem{proposition}[theorem]{Proposition}
\newtheorem{lemma}[theorem]{Lemma}
\newtheorem{conjecture}[theorem]{Conjecture}

\theoremstyle{remark}

\newtheorem{remark}[theorem]{Remark}

\theoremstyle{definition}
\newtheorem{definition}[theorem]{Definition}
\newtheorem{question}[theorem]{Question}

\begin{document}

\title{Unconditional Reflexive Polytopes} 

\author{Florian Kohl}
\address{Department of Mathematics and Systems Analysis\\
		Aalto University\\
		Espoo, Finland	}
\email{florian.kohl@aalto.fi}

\author{McCabe Olsen}
\address{Department of Mathematics\\
         The Ohio State University\\
         Columbus, OH 43210, USA}
\email{olsen.149@osu.edu}

\author{Raman Sanyal}
\address{Institut f{\"ur} Mathematik\\ 
	Goethe-Universit{\"a}t, Frankfurt am Main\\
	Germany}
\email{sanyal@math.uni-frankfurt.de}



\date{\today}

\thanks{
The first author was supported by Academy of
Finland project numbers 288318 and 13324921. The third author was supported by the
DFG-Collaborative Research Center, TRR 109 ``Discretization in Geometry and
Dynamics''.}


\begin{abstract}
    A convex body is unconditional if it is symmetric with respect to
    reflections in all coordinate hyperplanes. In this paper, we investigate
    unconditional lattice polytopes with respect to geometric, combinatorial,
    and algebraic properties. In particular, we characterize unconditional
    reflexive polytopes in terms of perfect graphs. As a prime example, we
    study the signed Birkhoff polytope. Moreover, we derive constructions for
    Gale-dual pairs of polytopes and we explicitly describe Gr\"obner bases
    for unconditional reflexive polytopes coming from partially ordered sets.
\end{abstract}

\maketitle


\section{Introduction}

\begin{figure}
\begin{center}
\begin{tikzpicture}[scale=0.6]
	\coordinate (a) at (-1,-1);
	\coordinate (b) at (-1,1);
	\coordinate (c) at (1,1);
	\coordinate (d) at (1,-1);
	\begin{pgftransparencygroup}
	\pgfsetfillopacity{0.2}
		\filldraw[line width=0.5mm, fill=gray] (a) -- (b) -- (c) -- (d) -- cycle;
	\end{pgftransparencygroup}
	\filldraw (0,0) circle (3pt);
   	\filldraw (-1,1) circle (3pt);
   	\filldraw (1,1) circle (3pt);
   	\filldraw (0,1) circle (3pt);
   	\filldraw (-1,0) circle (3pt);
   	\filldraw (1,0) circle (3pt);
   	\filldraw (-1,-1) circle (3pt);
   	\filldraw (0,-1) circle (3pt);
   	\filldraw (1,-1) circle (3pt);

	\draw [<->,>=stealth,thick] (2,0) -- (3,0);   	
   	
	\begin{pgftransparencygroup}
	\pgfsetfillopacity{0.2}
		\filldraw[line width=0.5mm, fill=gray] (4,0) -- (5,1) -- (6,0) -- (5,-1) -- cycle;
			\end{pgftransparencygroup}
	\filldraw (5,0) circle (3pt);   	
   	\filldraw (6,0) circle (3pt);
   	\filldraw (5,1) circle (3pt);
   	\filldraw (4,0) circle (3pt);
   	\filldraw (5,-1) circle (3pt);
\end{tikzpicture}\qquad
\begin{tikzpicture}[scale=0.6]
	\begin{pgftransparencygroup}
	\pgfsetfillopacity{0.2}
		\filldraw[line width=0.5mm, fill=gray] (0,1) -- (1,0) -- (-1,-1) -- cycle;
	\end{pgftransparencygroup}
	\filldraw (0,0) circle (3pt);
   	\filldraw (0,1) circle (3pt);
   	\filldraw (1,0) circle (3pt);
   	\filldraw (-1,-1) circle (3pt);

	\draw [<->,>=stealth,thick] (2,0) -- (3,0);

	\filldraw (5,0) circle (3pt);
    \begin{pgftransparencygroup}
	\pgfsetfillopacity{0.2}
		\filldraw[line width=0.5mm, fill=gray] (3,1) -- (6,1) -- (6,-2) -- cycle;
	\end{pgftransparencygroup}
   	\filldraw (3,1) circle (3pt);
   	\filldraw (4,1) circle (3pt);
   	\filldraw (5,1) circle (3pt);
   	\filldraw (6,1) circle (3pt);
   	\filldraw (4,0) circle (3pt);
   	\filldraw (6,0) circle (3pt);
   	\filldraw (5,-1) circle (3pt);
   	\filldraw (6,-1) circle (3pt);
   	\filldraw (6,-2) circle (3pt);
\end{tikzpicture}

\vspace{0.5cm}

\begin{tikzpicture}[scale=0.6]
	\begin{pgftransparencygroup}
	\pgfsetfillopacity{0.2}
		\filldraw[line width=0.5mm, fill=gray] (-2,-1) -- (0,1) -- (2,-1) -- cycle;
	\end{pgftransparencygroup}
	\filldraw (0,0) circle (3pt);
   	\filldraw (-2,-1) circle (3pt);
   	\filldraw (-1,-1) circle (3pt);
   	\filldraw (0,-1) circle (3pt);
   	\filldraw (1,-1) circle (3pt);
   	\filldraw (2,-1) circle (3pt);
   	\filldraw (1,0) circle (3pt);
   	\filldraw (0,1) circle (3pt);
   	\filldraw (-1,0) circle (3pt);

	\draw [<->,>=stealth,thick] (2,0) -- (3,0); 
	
	\begin{pgftransparencygroup}
	\pgfsetfillopacity{0.2}
		\filldraw[line width=0.5mm, fill=gray] (4,1) -- (6,1) -- (5,-1) -- cycle;
	\end{pgftransparencygroup} 
	\filldraw (5,0) circle (3pt);
   	\filldraw (4,1) circle (3pt);
   	\filldraw (5,1) circle (3pt);
   	\filldraw (6,1) circle (3pt);
   	\filldraw (5,-1) circle (3pt);
\end{tikzpicture}\qquad
\begin{tikzpicture}[scale=0.6]
	\begin{pgftransparencygroup}
	\pgfsetfillopacity{0.2}
		\filldraw[line width=0.5mm, fill=gray] (-1,-1) -- (2,-1) -- (0,1) -- (-1,1) -- cycle;
	\end{pgftransparencygroup} 
	\filldraw (0,0) circle (3pt);
   	\filldraw (-1,-1) circle (3pt);
   	\filldraw (0,-1) circle (3pt);
   	\filldraw (1,-1) circle (3pt);
   	\filldraw (2,-1) circle (3pt);
   	\filldraw (1,0) circle (3pt);
   	\filldraw (0,1) circle (3pt);
   	\filldraw (-1,1) circle (3pt);
   	\filldraw (-1,0) circle (3pt);
   	
   	\draw [<->,>=stealth,thick] (2,0) -- (3,0); 

	\filldraw (5,0) circle (3pt);
	\begin{pgftransparencygroup}
	\pgfsetfillopacity{0.2}
		\filldraw[line width=0.5mm, fill=gray] (4,0) -- (5,1) -- (6,1) -- (5,-1) -- cycle;
	\end{pgftransparencygroup}  
   	\filldraw (5,-1) circle (3pt);
   	\filldraw (6,1) circle (3pt);
   	\filldraw (5,1) circle (3pt);
   	\filldraw (4,0) circle (3pt);
\end{tikzpicture}

\vspace{0.5cm}

\begin{tikzpicture}[scale=0.6]
	\begin{pgftransparencygroup}
	\pgfsetfillopacity{0.2}
		\filldraw[line width=0.5mm, fill=gray] (-2,-1) -- (0,1) -- (1,0) -- (1,-1) -- cycle;
	\end{pgftransparencygroup} 
	\filldraw (0,0) circle (3pt);
   	\filldraw (-2,-1) circle (3pt);
   	\filldraw (-1,-1) circle (3pt);
   	\filldraw (0,-1) circle (3pt);
   	\filldraw (1,-1) circle (3pt);
   	\filldraw (1,0) circle (3pt);
   	\filldraw (0,1) circle (3pt);
   	\filldraw (-1,0) circle (3pt);
   	
   	\draw [<->,>=stealth,thick] (2,0) -- (3,0); 

	\begin{pgftransparencygroup}
	\pgfsetfillopacity{0.2}
		\filldraw[line width=0.5mm, fill=gray] (4,1) -- (6,1) -- (6,0) -- (5,-1)  -- cycle;
	\end{pgftransparencygroup}  
	\filldraw (5,0) circle (3pt);   	
   	\filldraw (5,-1) circle (3pt);
   	\filldraw (6,0) circle (3pt);
   	\filldraw (6,1) circle (3pt);
   	\filldraw (5,1) circle (3pt);
   	\filldraw (4,1) circle (3pt);
\end{tikzpicture}\qquad
\begin{tikzpicture}[scale=0.6]
	\begin{pgftransparencygroup}
	\pgfsetfillopacity{0.2}
		\filldraw[line width=0.5mm, fill=gray] (-1,1) -- (0,1) -- (1,0) -- (1,-1) -- (-1,-1) -- cycle;
	\end{pgftransparencygroup} 
	\filldraw (0,0) circle (3pt);
   	\filldraw (-1,-1) circle (3pt);
   	\filldraw (0,-1) circle (3pt);
   	\filldraw (1,-1) circle (3pt);
   	\filldraw (1,0) circle (3pt);
   	\filldraw (0,1) circle (3pt);
   	\filldraw (-1,1) circle (3pt);
   	\filldraw (-1,0) circle (3pt);
   	
   	\draw [<->,>=stealth,thick] (2,0) -- (3,0); 

	\begin{pgftransparencygroup}
	\pgfsetfillopacity{0.2}
		\filldraw[line width=0.5mm, fill=gray] (4,0) -- (5,1) -- (6,1) -- (6,0) -- (5,-1) -- cycle;
	\end{pgftransparencygroup}  
	\filldraw (5,0) circle (3pt);
   	\filldraw (5,-1) circle (3pt);
   	\filldraw (6,0) circle (3pt);
   	\filldraw (6,1) circle (3pt);
   	\filldraw (5,1) circle (3pt);
   	\filldraw (4,0) circle (3pt);
\end{tikzpicture}

\vspace{0.5cm}

\begin{tikzpicture}[scale=0.6]
\begin{pgftransparencygroup}
	\pgfsetfillopacity{0.2}
		\filldraw[line width=0.5mm, fill=gray] (-1,1) -- (0,1) -- (1,0) -- (1,-1) -- (0,-1) -- (-1,0) -- cycle;
	\end{pgftransparencygroup}  
	\filldraw (0,0) circle (3pt);
   	\filldraw (-1,1) circle (3pt);
   	\filldraw (0,1) circle (3pt);
   	\filldraw (1,0) circle (3pt);
   	\filldraw (1,-1) circle (3pt);
   	\filldraw (0,-1) circle (3pt);
   	\filldraw (-1,0) circle (3pt);
\end{tikzpicture}\qquad
\begin{tikzpicture}[scale=0.6]
\begin{pgftransparencygroup}
	\pgfsetfillopacity{0.2}
		\filldraw[line width=0.5mm, fill=gray] (-1,1) -- (0,1) -- (1,0) -- (0,-1) -- (-1,-1) -- cycle;
	\end{pgftransparencygroup}  
	\filldraw (0,0) circle (3pt);
   	\filldraw (-1,1) circle (3pt);
   	\filldraw (0,1) circle (3pt);
   	\filldraw (1,0) circle (3pt);
   	\filldraw (0,-1) circle (3pt);
   	\filldraw (-1,-1) circle (3pt);
   	\filldraw (-1,0) circle (3pt);
\end{tikzpicture}\qquad
\begin{tikzpicture}[scale=0.6]
\begin{pgftransparencygroup}
	\pgfsetfillopacity{0.2}
		\filldraw[line width=0.5mm, fill=gray] (-2,-1) -- (0,1) -- (1,0) -- (0,-1) -- cycle;
	\end{pgftransparencygroup}  
	\filldraw (0,0) circle (3pt);
   	\filldraw (-1,-1) circle (3pt);
   	\filldraw (0,-1) circle (3pt);
   	\filldraw (1,0) circle (3pt);
   	\filldraw (0,1) circle (3pt);
   	\filldraw (-1,0) circle (3pt);
   	\filldraw (-2,-1) circle (3pt);
\end{tikzpicture}\qquad
\begin{tikzpicture}[scale=0.6]
\begin{pgftransparencygroup}
	\pgfsetfillopacity{0.2}
		\filldraw[line width=0.5mm, fill=gray] (-2,-1) -- (0,1) -- (1,-1) -- cycle;
	\end{pgftransparencygroup}  
	\filldraw (0,0) circle (3pt);
   	\filldraw (-2,-1) circle (3pt);
   	\filldraw (-1,-1) circle (3pt);
   	\filldraw (0,-1) circle (3pt);
   	\filldraw (1,-1) circle (3pt);
   	\filldraw (0,1) circle (3pt);
   	\filldraw (-1,0) circle (3pt);
\end{tikzpicture}
\caption{All $16$ reflexive $2$-dimensional polytopes. This is Figure~1.5 in \cite{JanDiss}.}
\label{16reflexivePolygons-Picture}
\end{center}
\end{figure}
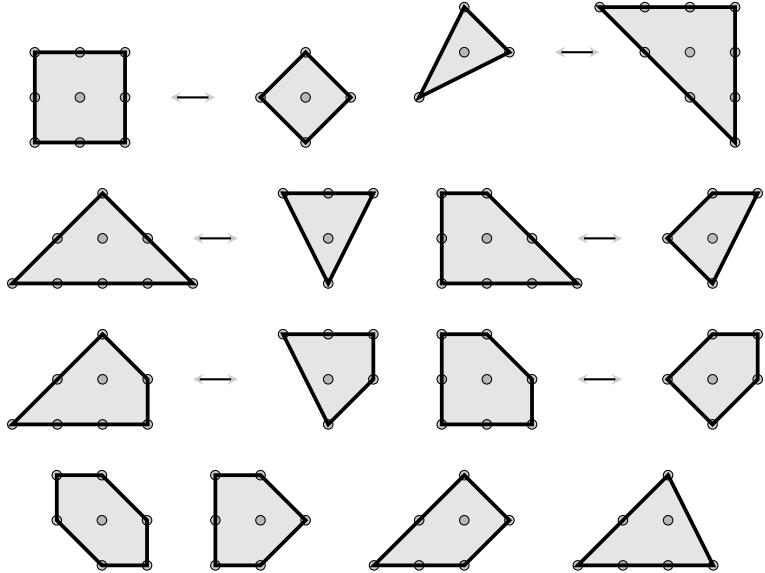

A $d$-dimensional convex lattice polytope $P \subset \R^d$ is called
\Def{reflexive} if its polar dual $P^*$ is again a lattice polytope.
Reflexive polytopes were introduced by Batyrev~\cite{Batyrev} in the context
of mirror symmetry as a reflexive polytope and its dual give rise to a
mirror-dual pair of Calabi--Yau manifolds (c.f.~\cite{CoxMirrorSurvey}).  As
thus, the results of Batyrev, and the subsequent connection with string
theory, have stimulated interest in the classification of reflexive polytopes
both among mathematical and theoretical physics communities.  As a consequence
of a well-known result of Lagarias and Ziegler~\cite{lagariasziegler}, there
are only finitely many reflexive polytopes in each dimension, up to unimodular
equivalence. In two dimensions, it is a straightforward exercise to verify
that there are precisely $16$ reflexive polygons, as depicted in
Figure~\ref{16reflexivePolygons-Picture}.  While still finite, there are
significantly more reflexive polytopes in higher dimensions.  Kreuzer and
Skarke~\cite{kreuzerskarke,kreuzerskarke4dim}~have completely classified
reflexive polytopes in dimensions $3$ and $4$, noting that there are exactly
$4319$ reflexive polytopes in dimension $3$ and  $473800776$ reflexive
polytopes in dimension $4$. The number of reflexive polytopes in dimension $5$
is not known.

In recent years, there has been significant progress in characterizing reflexive
polytopes in known classes of polytopes coming from combinatorics or
optimization; see, for example,~\cite{BeckHaaseSam09, Tagami10, Ohsugi14,
Hibi15-1, ChappellFriedlSanyal}. The purpose of this paper is to study a class
of reflexive polytopes motivated by convex geometry and relate it to
combinatorics. A convex body $K \subset \R^d$ is \Def{unconditional} if $\p
\in K$ if and only if $\sigma \p := (\sigma_1 p_1,\sigma_2 p_2, \dots,\sigma_d
p_d) \in K$ for all $\sigma \in \PM^d$. Unconditional convex bodies, for
example, arise as unit balls in the theory of Banach spaces with a
$1$-unconditional basis. They constitute a restricted yet surprisingly
interesting class of convex bodies for which a number of claims have been
verified; cf.~\cite{BGVV}. For example, we mention that the \emph{Mahler
conjecture} is known to hold for unconditional convex bodies; see
Section~\ref{sec:OLPUnconditional}. In this paper, we investigate
unconditional lattice polytopes and their relation to anti-blocking polytopes
from combinatorial optimization. In particular, we completely characterize
unconditional reflexive polytopes.\\

The structure of this paper is as follows. In Section~\ref{sec:Background}, we
briefly review notions and results from discrete geometry and Ehrhart theory.

In Section~\ref{sec:OLPUnconditional}, we introduce and study unconditional
and, more generally, locally anti-blocking polytopes. The main result is
Theorem~\ref{thm:ExtendTriang} that relates regular, unimodular, and flag
triangulations to the associated anti-blocking polytopes. 

In Section~\ref{sec:graphs}, we associate an unconditional lattice polytope
$\U P_G$ to every finite graph $G$. We show in
Theorems~\ref{thm:unconditional} and~\ref{thm:uncond_refl_different} that an
unconditional polytope $P$ is reflexive if and only if $P = \U P_G$ for some
unique perfect graph $G$. This also implies that unconditional reflexive
polytopes have regular, unimodular triangulations.

Section~\ref{sec:BBn} is devoted to a particular family of unconditional
reflexive polytopes and is of independent interest: We show that the
\Def{type-B Birkhoff polytope} or \Def{signed Birkhoff polytope} $\BB(n)$,
that is, the convex hull of signed permutation matrices, is an unconditional
reflexive polytope.  We compute normalized volumes and $h^*$-vectors of
$\BB(n)$ and its dual $\Cn \coloneqq \BB(n)^*$ for small values of $n$. 

The usual Birkhoff polytope and the Gardner polytope of~\cite{GardnerPolytope}
appear as faces of $\BB(n)$ and $\Cn$, respectively. These two polytopes form
a \emph{Gale-dual pair} in the sense of~\cite{GardnerPolytope}. 
In Section~\ref{sec:gale}, we give a general construction for compressed
Gale-dual pairs coming from CIS graphs.

In Section~\ref{sec:ChainAndGroebner}, we investigate unconditional polytopes
associated to comparability graphs of posets. In particular, we explicitly
describe a quadratic square-free Gr\"obner basis for the corresponding toric
ideal. 

We close with open questions and future directions in
Section~\ref{sec:conclude}. \\

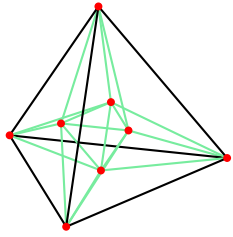
\begin{figure}

\begin{tikzpicture}[x  = {(0.963586094867353cm,-0.100597083561054cm)},
                    y  = {(0.142869993056962cm,0.976884553520567cm)},
                    z  = {(-0.226030977660555cm,0.188618121801521cm)},
                    scale = 1.5,
                    color = {lightgray}]

  \definecolor{pointcolor_p}{rgb}{ 1,0,0 }
  \tikzstyle{pointstyle_p} = [fill=pointcolor_p]

  \coordinate (v0_p) at (0.310345, 0, 0.219447);
  \coordinate (v1_p) at (-0.310345, 0, 0.219447);
  \coordinate (v2_p) at (1, 0, -0.707107);
  \coordinate (v3_p) at (-1, 0, -0.707107);
  \coordinate (v4_p) at (0, 0.310345, -0.219447);
  \coordinate (v5_p) at (0, -0.310345, -0.219447);
  \coordinate (v6_p) at (0, 1, 0.707107);
  \coordinate (v7_p) at (0, -1, 0.707107);

  \definecolor{linecolor_p_v1_v0}{rgb}{ 0.4667,0.9255,0.6196 }
  \definecolor{linecolor_p_v2_v0}{rgb}{ 0.4667,0.9255,0.6196 }
  \definecolor{linecolor_p_v3_v1}{rgb}{ 0.4667,0.9255,0.6196 }
  \definecolor{linecolor_p_v3_v2}{rgb}{ 0,0,0 }
  \definecolor{linecolor_p_v4_v0}{rgb}{ 0.4667,0.9255,0.6196 }
  \definecolor{linecolor_p_v4_v1}{rgb}{ 0.4667,0.9255,0.6196 }
  \definecolor{linecolor_p_v4_v2}{rgb}{ 0.4667,0.9255,0.6196 }
  \definecolor{linecolor_p_v4_v3}{rgb}{ 0.4667,0.9255,0.6196 }
  \definecolor{linecolor_p_v5_v0}{rgb}{ 0.4667,0.9255,0.6196 }
  \definecolor{linecolor_p_v5_v1}{rgb}{ 0.4667,0.9255,0.6196 }
  \definecolor{linecolor_p_v5_v2}{rgb}{ 0.4667,0.9255,0.6196 }
  \definecolor{linecolor_p_v5_v3}{rgb}{ 0.4667,0.9255,0.6196 }
  \definecolor{linecolor_p_v5_v4}{rgb}{ 0.4667,0.9255,0.6196 }
  \definecolor{linecolor_p_v6_v0}{rgb}{ 0.4667,0.9255,0.6196 }
  \definecolor{linecolor_p_v6_v1}{rgb}{ 0.4667,0.9255,0.6196 }
  \definecolor{linecolor_p_v6_v2}{rgb}{ 0,0,0 }
  \definecolor{linecolor_p_v6_v3}{rgb}{ 0,0,0 }
  \definecolor{linecolor_p_v6_v4}{rgb}{ 0.4667,0.9255,0.6196 }
  \definecolor{linecolor_p_v7_v0}{rgb}{ 0.4667,0.9255,0.6196 }
  \definecolor{linecolor_p_v7_v1}{rgb}{ 0.4667,0.9255,0.6196 }
  \definecolor{linecolor_p_v7_v2}{rgb}{ 0,0,0 }
  \definecolor{linecolor_p_v7_v3}{rgb}{ 0,0,0 }
  \definecolor{linecolor_p_v7_v5}{rgb}{ 0.4667,0.9255,0.6196 }
  \definecolor{linecolor_p_v7_v6}{rgb}{ 0,0,0 }

  \tikzstyle{linestyle_p_v1_v0} = [color=linecolor_p_v1_v0, thick]
  \tikzstyle{linestyle_p_v2_v0} = [color=linecolor_p_v2_v0, thick]
  \tikzstyle{linestyle_p_v3_v1} = [color=linecolor_p_v3_v1, thick]
  \tikzstyle{linestyle_p_v3_v2} = [color=linecolor_p_v3_v2, thick]
  \tikzstyle{linestyle_p_v4_v0} = [color=linecolor_p_v4_v0, thick]
  \tikzstyle{linestyle_p_v4_v1} = [color=linecolor_p_v4_v1, thick]
  \tikzstyle{linestyle_p_v4_v2} = [color=linecolor_p_v4_v2, thick]
  \tikzstyle{linestyle_p_v4_v3} = [color=linecolor_p_v4_v3, thick]
  \tikzstyle{linestyle_p_v5_v0} = [color=linecolor_p_v5_v0, thick]
  \tikzstyle{linestyle_p_v5_v1} = [color=linecolor_p_v5_v1, thick]
  \tikzstyle{linestyle_p_v5_v2} = [color=linecolor_p_v5_v2, thick]
  \tikzstyle{linestyle_p_v5_v3} = [color=linecolor_p_v5_v3, thick]
  \tikzstyle{linestyle_p_v5_v4} = [color=linecolor_p_v5_v4, thick]
  \tikzstyle{linestyle_p_v6_v0} = [color=linecolor_p_v6_v0, thick]
  \tikzstyle{linestyle_p_v6_v1} = [color=linecolor_p_v6_v1, thick]
  \tikzstyle{linestyle_p_v6_v2} = [color=linecolor_p_v6_v2, thick]
  \tikzstyle{linestyle_p_v6_v3} = [color=linecolor_p_v6_v3, thick]
  \tikzstyle{linestyle_p_v6_v4} = [color=linecolor_p_v6_v4, thick]
  \tikzstyle{linestyle_p_v7_v0} = [color=linecolor_p_v7_v0, thick]
  \tikzstyle{linestyle_p_v7_v1} = [color=linecolor_p_v7_v1, thick]
  \tikzstyle{linestyle_p_v7_v2} = [color=linecolor_p_v7_v2, thick]
  \tikzstyle{linestyle_p_v7_v3} = [color=linecolor_p_v7_v3, thick]
  \tikzstyle{linestyle_p_v7_v5} = [color=linecolor_p_v7_v5, thick]
  \tikzstyle{linestyle_p_v7_v6} = [color=linecolor_p_v7_v6, thick]

  \draw[linestyle_p_v1_v0] (v1_p) -- (v0_p);
  \draw[linestyle_p_v2_v0] (v2_p) -- (v0_p);
  \draw[linestyle_p_v3_v1] (v3_p) -- (v1_p);
  \draw[linestyle_p_v3_v2] (v3_p) -- (v2_p);
  \draw[linestyle_p_v4_v0] (v4_p) -- (v0_p);
  \draw[linestyle_p_v4_v1] (v4_p) -- (v1_p);
  \draw[linestyle_p_v4_v2] (v4_p) -- (v2_p);
  \draw[linestyle_p_v4_v3] (v4_p) -- (v3_p);
  \draw[linestyle_p_v5_v0] (v5_p) -- (v0_p);
  \draw[linestyle_p_v5_v1] (v5_p) -- (v1_p);
  \draw[linestyle_p_v5_v2] (v5_p) -- (v2_p);
  \draw[linestyle_p_v5_v3] (v5_p) -- (v3_p);
  \draw[linestyle_p_v5_v4] (v5_p) -- (v4_p);
  \draw[linestyle_p_v6_v0] (v6_p) -- (v0_p);
  \draw[linestyle_p_v6_v1] (v6_p) -- (v1_p);
  \draw[linestyle_p_v6_v2] (v6_p) -- (v2_p);
  \draw[linestyle_p_v6_v3] (v6_p) -- (v3_p);
  \draw[linestyle_p_v6_v4] (v6_p) -- (v4_p);

  \fill[pointcolor_p] (v4_p) circle (1 pt);
  \node at (v4_p) [text=black, inner sep=0.5pt, above right, draw=none, align=left] {};

  \draw[linestyle_p_v7_v0] (v7_p) -- (v0_p);

  \fill[pointcolor_p] (v0_p) circle (1 pt);
  \node at (v0_p) [text=black, inner sep=0.5pt, above right, draw=none, align=left] {}; 

  \draw[linestyle_p_v7_v1] (v7_p) -- (v1_p);

  \fill[pointcolor_p] (v1_p) circle (1 pt);
  \node at (v1_p) [text=black, inner sep=0.5pt, above right, draw=none, align=left] {}; 

  \draw[linestyle_p_v7_v2] (v7_p) -- (v2_p);

  \fill[pointcolor_p] (v2_p) circle (1 pt);
  \node at (v2_p) [text=black, inner sep=0.5pt, above right, draw=none, align=left] {};

  \draw[linestyle_p_v7_v3] (v7_p) -- (v3_p);

  \fill[pointcolor_p] (v3_p) circle (1 pt);
  \node at (v3_p) [text=black, inner sep=0.5pt, above right, draw=none, align=left] {}; 

  \draw[linestyle_p_v7_v5] (v7_p) -- (v5_p);

  \fill[pointcolor_p] (v5_p) circle (1 pt);
  \node at (v5_p) [text=black, inner sep=0.5pt, above right, draw=none, align=left] {}; 

  \draw[linestyle_p_v7_v6] (v7_p) -- (v6_p);

  \fill[pointcolor_p] (v7_p) circle (1 pt);
  \node at (v7_p) [text=black, inner sep=0.5pt, above right, draw=none, align=left] {}; 
  \fill[pointcolor_p] (v6_p) circle (1 pt);
  \node at (v6_p) [text=black, inner sep=0.5pt, above right, draw=none, align=left] {}; 


\end{tikzpicture}
\caption{Schlegel diagram for $\BB(2)$.}
\label{fig:BB2}
\end{figure}

\noindent
\textbf{Acknowledgements.} The first two authors would like to thank Matthias
Beck, Benjamin Braun, and Jan Hofmann for helpful comments and suggestions for
this work. Furthermore, Figure~\ref{fig:BB2} was created by Benjamin
Schr\"oter. Additionally, the authors thank Takayuki Hibi and Akiyoshi
Tsuchiya for organizing the 2018 Summer Workshop on Lattice Polytopes at Osaka
University where this work began. The third author thanks Kolja Knauer and
Sebastian Manecke for insightful conversations.

\section{Background}
\label{sec:Background}

In this section, we provide a brief introduction to polytopes and Ehrhart
theory. For additional background and details, we refer the reader to the
excellent books \cite{CCD,Ziegler}. A \Def{polytope} in $\R^d$ is the
inclusion-minimal convex set $P = \conv(\vb_1,\dots,\vb_n)$ containing a given
collection of points $\vb_1,\dots,\vb_n \in \R^d$. The unique
inclusion-minimal set $V \subseteq P$ such that $P = \conv(V)$ is called the
\Def{vertex set} and is denoted by $V(P)$. If $V(P) \subset
\Z^d$, then $P$ is called a \Def{lattice polytope}. By the Minkowski--Weyl theorem,
polytopes are precisely the bounded sets of the form
\[
    P \ = \ \{ \x \in \R^d : \inner{\a_i,\x} \le b_i \text{ for } i = 1,\dots,m \}
\]
for some $\a_1,\dots,\a_m \in \R^d$ and $b_1,\dots,b_m \in \R$. If
$\inner{\a_i,\x} \le b_i$ is irredundant, then $F = P \cap \{\x : \inner{\a_i,\x} =
b_i \}$ is a facet and the inequality is said to be \Def{facet-defining}.

The \Def{dimension} of a polytope $P$ is defined to be the dimension of its
affine span. A $d$-dimensional polytope has at least $d+1$ vertices and a
$d$-polytope with exactly $d+1$ many vertices is called a \Def{$d$-simplex}.
A $d$-simplex $\Delta = \conv\{\vb_0, \vb_1, \dots, \vb_d\}$ is called
\Def{unimodular} if $\vb_1 - \vb_0$, $\vb_2 - \vb_0$, $\dots$, $\vb_d - \vb_0$
form a basis for the lattice $\Z^d$, or equivalently if
$\vol(\Delta)=\frac{1}{d!}$, where $\vol$ is the Euclidean volume. For lattice
polytopes $P \subset \R^d$, we define the \Def{normalized volume} $\Vol(P)
\coloneqq d! \vol(P)$. So unimodular simplices are the lattice polytopes with
normalized volume $1$.  We say that two lattice polytopes $P,P' \subset \R^d$
are \Def{unimodularly equivalent} if $P' = T(P)$  for some transformation
$T(\x) = W\x + \vb$ with $W \in \operatorname{SL}_d(\Z)$ and $\vb\in
\Z^d$.	In particular, any two unimodular simplices are unimodularly
equivalent.

Given a lattice $e$-polytope $P \subset \R^d$ and $t\in\Zpos$, let $t P
\coloneqq \{t\cdot \x \, : \, \x \in P\}$ be the $t^{\text{th}}$ dilate of
$P$.  By a famous result of Ehrhart \cite[Thm.~1]{Ehrhart}, the lattice-point
enumeration function 
\[
    \ehr_{P}(t) \ \coloneqq \ |t P \cap \Z^d|
\]
agrees with a polynomial in the variable $t$ of degree $e$ with leading
coefficient $\vol(P)$ and is called the \Def{Ehrhart polynomial}.   This also
implies that the formal generating function
\[
    1 + \sum_{t \geq 1} \ehr_{P}(t) z^t \ = \ \frac{h^{\ast}_0(P) +
    h^{\ast}_1(P) z + \dots + h^{\ast}_e(P) z^e }{(1-z)^{e+1}}
\]
is a rational function such that the degree of the numerator is at most $e$
(see, e.g., {\cite[Lem 3.9]{CCD}}). We call the numerator the
\Def{$h^{\ast}$-polynomial} of $P$. The vector $h^{\ast}(P) =
(h_0^\ast(P),h_1^\ast(P),\ldots,h_d^\ast(P)) \in \Z^{d+1}$, where we set
$h_i^\ast(P) := 0$ for $\dim P < i \le d$, is called the
\Def{$h^{\ast}$-vector} of $P$.  One should note that the Ehrhart polynomial
is invariant under unimodular transformations.  

\begin{theorem}[{\cite{StanleyNonnegativity, Stan93}}]\label{thm:h-mono}
    Let $P \subseteq Q \subset \R^d$ be lattice polytopes. Then 
    \[
        0 \ \le \ h^*_i(P) \ \le \ h^*_i(Q)
    \]
    for all $i = 0,\dots,d$.
\end{theorem}

The $h^{\ast}$-vector encodes significant information about the underlying
polytope. This is nicely illustrated in the case of reflexive polytopes.  For
a $d$-polytope $P \subset \R^d$ with $0$ in the interior, we define the
\Def{(polar) dual polytope}
\[
    P^\ast \ \coloneqq \ \{ \y \in \R^d \ : \ \inner{\y, \x} \leq
    1 \ \text{for all } \x \in P \} \, .
\]

\begin{definition}\label{def:reflexive}
    Let $P\subset \R^d$ be a $d$-dimensional lattice polytope that contains
    the origin in its interior. We say that $P$ is \Def{reflexive} if $P^*$ is
    also a lattice polytope. Equivalently, $P$ is reflexive if it has a
    description of the form
    \[
    P  \ = \ \{ \x \in \R^d \colon  \inner{\a_i,\x} \leq 1 \text{ for } i =
    1,\dots,m \},
    \]
    for some $\a_1,\dots,\a_m \in \Z^d$.
\end{definition}

Reflexivity can be completely characterized by enumerative data  of the
$h^{\ast}$-vector.
\begin{theorem}[{\cite[Thm.~2.1]{HibiDualPolytopes}}]
    Let $P \subset \R^d$ be a $d$-dimensional lattice polytope with $h^*(P) =
    (h_0^*,\dots,h_d^*)$.  Then $P$ is unimodularly equivalent to a reflexive
    polytope if and only if $h^{\ast}_k = h^{\ast}_{d-k}$ for all $0 \leq k
    \leq \left \lfloor \frac d 2 \right \rfloor$.
\end{theorem}

The reflexivity property is also deeply related to commutative algebra. A
polytope $P \subset \R^d$ is reflexive if the canonical module of the
associated graded algebra $\k[P]$ is (up to a shift in grading) isomorphic to
$\k[P]$ and its minimal generator has degree $1$. If one allows the unique
minimal generator to have arbitrary degree, one arrives at the notion of
Gorenstein rings, for details we refer to \cite[Sec 6.C]{BrunsGubeladze}.  We
say that $P$ is \Def{Gorenstein} if there exist $c\in\Zpos$ and $\q \in
\Z^d$ such that $\q + c P$ is a reflexive polytope. This is equivalent to saying
that $\k[P]$ is Gorenstein. The dilation factor $c$ is often called
the \Def{codegree}. In particular, reflexive polytopes are Gorenstein of
codegree $1$. By combining results of Stanley \cite{StanleyGorenstein} and De
Negri--Hibi \cite{DeNegriHibi}, we have a characterization of the Gorenstein
property in terms of the $h^\ast$-vector. Namely, $P$ is Gorenstein if
and only if $h^\ast_i=h^\ast_{d-c+1-i}$ for all $i$.

Aside from examining algebraic properties of lattice polytopes, one can also
investigate discrete geometric properties.  Every lattice polytope admits a
subdivision into lattice simplices.  Even more, one can guarantee that every
lattice point contained in a polytope corresponds to a vertex of such a
subdivision.  However, one cannot guarantee the existence of a subdivision
where all simplices are unimodular when the dimension is greater than $2$.
This leads us to our next definition:

\begin{definition} \label{def:ch1UnimodularTriangulation}
    Let  $P$ be a $d$-dimensional lattice polytope given by $P = \conv(V)$ for
    some finite set $V \subset \Z^d$. A \Def{subdivision} of $P$ with vertices
    in $V$ is a collection $\Subdiv = \{P_1,\dots,P_r\}$ of $d$-dimensional
    polytopes with vertices in $V$ such that $P = P_1 \cup \cdots \cup P_r$
    and $P_i \cap P_j$ is a common face of $P_i$ and $P_j$ for all $1 \le i <
    j \le r$. If all polytopes $P_i$ are (unimodular) simplices, then $\Subdiv$
    is a \Def{(unimodular) triangulation}. A subdivision $\Subdiv' =
    \{P_1',\dots,P_s'\}$ \Def{refines} $\Subdiv$ if for every $P_i' \in
    \Subdiv'$ there is $P_j \in \Subdiv$ such that $P_i' \subseteq P_j$.

\end{definition}

Suppose that $P = \conv(V)$. Any map $\omega : V \to \R$ yields a
piecewise-linear and convex function $\OV\omega : P \to \R$ by $\OV\omega(p)
:= \min\{ \lambda : (p,\lambda) \in P^\omega \}$, where $P^\omega = \conv(
(v,\omega(v)) : v \in V)$. The domains of linearity of $\OV\omega$ determine a
subdivision of $P$, called a \Def{regular} subdivision. The function $\omega$
is referred to as \Def{heights} on $P$. For more on (regular) subdivisions and
details we refer to~\cite{deloerarambausantos,HPPS}.

A particularly simple example of a regular triangulation of a (lattice) polytope $P$ is
the \Def{pulling triangulation}. For an arbitrary but fixed ordering of the
vertices $V(P)$, let $\vb_F$ be the first vertex in the face $F$ in the given
ordering. The pulling triangulation of $P$ is then defined recursively as
follows: If $P$ is a simplex, then $\Triang = \{ P \}$. Otherwise, let
$F_1,\dots,F_r$ be the facets of $P$ not containing $\vb_P$ and let
$\Triang_1,\dots,\Triang_r$ be their pulling triangulations with
respect to the induced ordering. Then 
\[
    \Triang \ := \ \{ \conv(\{\vb_P\} \cup S) : S \in \Triang_j, j=1,\dots,r \}
\]
is the pulling triangulation of $P$; see~\cite[Ch.~5.7]{BeckSanyal}.

For a subdivision $\Subdiv$, we write $V(\Subdiv) = \bigcup_{Q \in \Subdiv}
V(Q)$. If $\Triang$ is a triangulation, then $W \subseteq V(\Triang)$ is
called a \Def{non-face} if $\conv(W)$ is not a face of any $Q \in
\Triang$.  The triangulation $\Triang$ is called \Def{flag} if the
inclusion-minimal non-faces $W$ satisfy $|W| = 2$.

A special class of polytopes which possess regular, unimodular
triangulations are {compressed polytopes}. A polytope $P$ is
\Def{compressed} if every pulling triangulation is unimodular
\cite{StanleyNonnegativity}.  In the interest of providing a useful
characterization of compressed polytopes, we must define the notion of width
with respect to a facet. Let $P \subset \R^d$ be a $d$-dimensional lattice
polytope and $F_i = P \cap \{ \x : \inner{\a_i,\x} = b_i \}$ a facet. We
assume that $\a_i$ is primitive, that is, its coordinates are coprime. The
\Def{width} of $P$ with respect to the facet $F_i$ is
\[
    \max_{\p \in P} \inner{\a_i,\p}   - \min_{\p \in P} \inner{\a_i,\p}   \, .
\]
The maximum over all facets is called the \Def{facet width} of $P$.
\begin{theorem}[{\cite[Thm. 1.1]{OhsugiHibi-RevLex}} {\cite[Thm.
    2.4]{sullivant}}] \label{thm:compressed}
    Let $P \subset \R^d$ be a full-dimensional lattice polytope. The
    following are equivalent:
    \begin{enumerate}
		\item $P$ is compressed;
        \item $P$ has facet width one; 
        \item $P$ is unimodularly equivalent to the intersection of a unit
            cube with an affine space.
    \end{enumerate}
\end{theorem}

\begin{definition}
    A lattice polytope $P \subset \R^d$ has the \Def{integer decomposition
    property} (IDP) if for any positive integer $t$ and for all $\x\in
    tP\cap \Z^d$, there exist $\vb_1,\ldots,\vb_t \in P \cap \Z^d$
    such that $\x=\vb_1+\cdots+\vb_t$.
\end{definition}

One should note that if $P$ has a unimodular triangulation, then
$P$ has the IDP. However, there are examples of polytopes which have
the IDP, yet do not even admit a \Def{unimodular cover}, that is, a covering
of $P$ by unimodular simplices, see \cite[Sec. 3]{BrunsGubeladzeCounterexample}. A more complete hierarchy of covering
properties can be found in \cite{HPPS}.

We say that $h^\ast(P)$ is \Def{unimodal} if there exists a $k$ such that
$h^\ast_0\leq h^\ast_1\leq \cdots \leq h^\ast_k \geq \cdots \geq
h^\ast_{d-1}\geq h^\ast_d$. Unimodality appears frequently in combinatorial
settings and it often hints at a deeper underlying algebraic structure, see
\cite{AHK,Brenti-Unimodal,LogConcaveStanley}. One famous instance is given by
Gorenstein polytopes that admit a regular, unimodular triangulation. 

\begin{theorem}[{\cite[Thm. 1]{BrunsRoemer}}] \label{thm:BrunsRoemer}
    If $P$ is Gorenstein and has a regular, unimodular triangulation,
    then $h^\ast(P)$ is unimodal.
\end{theorem}

The following conjecture is commonly attributed to Ohsugi and Hibi \cite{OhsugiHibi-Conjecture}:
\begin{conjecture} \label{conj:OhsugiHibi}
    If $P$ is Gorenstein and has the IDP, then $h^\ast(P)$ is unimodal. 
\end{conjecture}

\section{Unconditional and anti-blocking polytopes}
\label{sec:OLPUnconditional}

For $\sigma \in \PM^d$ and $\p \in \R^d$, let us write $\sigma\p \ceq
(\sigma_1p_1,\sigma_2p_2,\dots,\sigma_dp_d)$. A convex polytope $P \subset
\R^d$ is called \Def{$1$-unconditional} or simply \Def{unconditional} if $\p
\in P$ implies $\sigma \p \in P$ for all $\sigma \in \PM^d$. So, unconditional
polytopes are precisely the polytopes that are invariant under reflection in
all coordinate hyperplanes.  It is apparent that $P$ can be recovered from its
restriction to the first orthant $\R^d_+ \ceq \{ x \in \R^d : x_1,\dots,x_d \ge
0\}$, which we denote by $P_+ \ceq P \cap \R^d_+$. The polytope $P_+$ has the
property that for any $\q\in P_+$  and $\p \in \R^d$ with $0 \le p_i \le q_i$
for all $i$, it holds that $\p \in P_+$. Polytopes in $\R^d_+$ with this
property are called \Def{anti-blocking polytopes}. Anti-blocking polytopes
were studied and named by Fulkerson~\cite{Fulkerson1,Fulkerson2} in the
context of combinatorial optimization, but they are also known as \emph{convex
corners} or \emph{down-closed} polytopes; see, for
example,~\cite{BollobasBrightwell}. 

Let us also write $\OV\p \ceq (|p_1|,|p_2|,\dots,|p_d|)$.  Given an anti-blocking
polytope $Q \subset \R^d_+$ it is straightforward to verify that 
\[
    \U Q \ \coloneqq \ \{ \p \in \R^d : \OV\p \in Q \}
\]
is an unconditional convex body.  Following Schrijver's treatment of
anti-blocking polytopes in~\cite[Sec.~9.3]{schrijver1986}, we recall that
every full-dimensional anti-blocking polytope has an irredundant inequality
description of the form
\begin{equation}\label{eqn:AB-ineqs}
    Q \ = \ \{ \x \in \R^d_+ : \inner{\a_i, \x} \le 1 \text{ for } i =
    1,\dots,m \}
\end{equation}
for some $\a_1,\dots,\a_m \in \R^d_+$. Also, we define 
\[
    \convDown{\c_1,\dots,\c_r} \ \coloneqq \ \R^d_+ \cap (
    \conv(\c_1,\dots,\c_r) + (-\R^d_+) ),
\]
where '$+$' denotes vector sum,
as the inclusion-minimal anti-blocking polytope containing the points $\c_1,\dots,\c_r \in \R^d_+$. 
Conversely, if we define $\Vdown(Q) := \{\vb_1,\dots,\vb_r\} $ to be the
vertices of an anti-blocking polytope $Q$ that are maximal with respect to the componentwise order, then
$Q = \convDown{\vb_1,\dots,\vb_r}$. We record the consequences for the
unconditional polytopes.

\begin{proposition}\label{prop:AB-inequalities}
    Let $P \subset \R^d_+$ be an anti-blocking $d$-polytope given
    by~\eqref{eqn:AB-ineqs}. Then an irredundant inequality description of $\U
    P$ is given by the distinct
    \[
        \inner{\sigma \a_i, \x } \le 1
    \]
    for $i=1,\dots,m$   and $\sigma \in \PM^d$. Likewise, the vertices of $\U
    P$ are $V(\U P) = \{ \sigma\vb : \vb \in \Vdown(P), \sigma \in \PM^d \}$.
\end{proposition}

Our first result relates properties of subdivisions of anti-blocking polytopes
to that of the associated unconditional polytopes.  The $2^d$ orthants in
$\R^d$ are denoted by $\R^d_\sigma \coloneqq \sigma \R^d_+$ for $\sigma \in \PM^d$.

\begin{theorem} \label{thm:ExtendTriang}
    Let $P \subset \R^d_+$ be an anti-blocking polytope with triangulation
    $\Triang$. Then 
    \[
        \U \Triang \ \coloneqq \ \{ \sigma S : S \in \Triang, \sigma \in \PM^d \}
    \]
    is a triangulation of $\U P$. Furthermore
    \begin{enumerate}[\rm (i)]
        \item If $\Triang$ is unimodular, then so is $\U\Triang$.
        \item If $\Triang$ is regular, then so is $\U\Triang$.
        \item If $\Triang$ is flag, then so is $\U\Triang$.
    \end{enumerate}
\end{theorem}

\begin{proof}
    It is clear that $\U \Triang$ is a triangulation of $\U P$ and 
    statement (i) is obvious. 
    
    To show (iii), let us assume that $\Triang$ is flag and let $W \subseteq
    V(\U\Triang)$ be an inclusion-minimal non-face of $\U\Triang$. If $W
    \subset \R^d_\sigma$ for some $\sigma$, then $\sigma W \subseteq
    V(\Triang)$ is an inclusion-minimal non-face of $\Triang$ and hence $|W| =
    |\sigma W| = 2$. Thus, there is $1 \le i \le d$ and $\p,\p' \in W$ with
    $p_i < 0 < p'_i$. But then $W' = \{ \p,\p' \} \subseteq W$ is also a
    non-face. Since we assume $W$ to be inclusion-minimal, this proves (iii).
    
    To show (ii), assume that $\Triang$ is regular and let $\omega : V(P) \to
    \R$ the corresponding heights. We extend $\omega$ to $V \coloneqq
    \bigcup_\sigma \sigma V(\Triang)$ by setting $\omega'(\vb) \coloneqq
    \|\vb\|_1 + \epsilon \, \omega(\OV\vb)$, where $\|\vb\|_1 = \sum_i |v_i|$
    and $\OV\vb = (|v_1|, \dots, |v_d|)$.  For $\epsilon = 0$ it is easy to
    see that the heights induce a regular subdivision of $P$ into $\sigma P$
    for $\sigma \in \PM^d$. For $\epsilon > 0$ sufficiently small, the heights
    $\omega'$ then induce the triangulation $\sigma \Triang$ on $\sigma P$.
\end{proof}

We call a polytope $P \subset \R^d$ \Def{locally anti-blocking} if $(\sigma P)
\cap \R^d_+$ is an  anti-blocking polytope for every $\sigma \in \PM^d$. In
particular, every locally anti-blocking polytope comes with a canonical
subdivision into polytopes $P_\sigma := P \cap \R^d_\sigma$ for $\sigma \in
\PM^d$.  Unconditional polytopes as well as anti-blocking polytopes are
clearly locally anti-blocking. It follows
from~\cite[Lemma~3.12]{ChappellFriedlSanyal} that for any two anti-blocking
polytopes $P_1, P_2 \subset \R^d_+$, the polytopes
\[
    P_1 + (-P_2) \qquad \text{ and }  \qquad P_1 \vee (-P_2) \ := \ \conv( P_1
    \cup -P_2 )
\]
are locally anti-blocking. Locally anti-blocking polytopes are studied in
depth in~\cite{LocallyAntiblocking}. The following is a simple, but important
observation.

\begin{lemma}\label{lem:LAB-reflexive}
    Let $P \subset \R^d$ be a locally anti-blocking lattice polytope with $0$
    in the interior. Then $P$ is reflexive if and only if $P_\sigma := P \cap
    \R^d_\sigma$ is compressed for all $\sigma \in \PM^d$.
\end{lemma}
\begin{proof}
    Since $P$ is a lattice polytope with $0$ in the interior, we can assume
    that $P$ is given as
    \[
        P \ = \ \{ \x \in \R^d : \inner{\a_i,\x} \le b_i \text{ for }
        i=1,\dots,m\}
    \]
    for some $\a_1,\dots,\a_m \in \Z^d$ primitive and $b_1,\dots,b_m \in
    \Q_{>0}$. Note that for $\sigma \in \PM^d$ we have that $P_\sigma  =  P
    \cap \R^d_\sigma$ is given by all $\x \in \R^d_\sigma$ such that 
    \[
        \inner{\a_i,\x} \le b_i \text{ for } i \in I_\sigma \, ,
    \]
    where $I_\sigma \subseteq \{ i \in [m] : \a_i \in \R^d_\sigma \}$ so that
    the inequalities are facet-defining. For $\sigma \in \PM^d$ and
    $i \in I_\sigma$, we observe 
    \[
        \min_{\p \in P_\sigma} \inner{\a_i,\p} \ = \ 0 \quad \text{ and }
        \quad
        \max_{\p \in P_\sigma} \inner{\a_i,\p} \ = \ b_i \, .
    \]
    The former holds because $\inner{\a_i,\p} \ge 0$ for $\a_i,\p \in
    \R^d_\sigma$ and the latter because $\inner{\a_i,\p} = b_i$ is
    facet-defining.  Hence $P_\sigma$ has facet width one with respect to
    $\a_i$ if and only if $b_i = 1$.  Since for every $i \in [m]$ there is a
    $\sigma \in \PM^d$ with $i \in I_\sigma$, we get using
    Theorem~\ref{thm:compressed} that if $P_\sigma$ is compressed, then $P$
    has a representation as in Definition~\ref{def:reflexive} and is therefore
    reflexive.

    Now let $P$ be reflexive and $\sigma$ arbitrary. We only need to argue
    that the facets of $P_\sigma$ corresponding to $\{ \x : x_j = 0 \}$ have
    facet width one. Assume that there is a vertex $\vb \in P_\sigma$ with
    $\sigma_jv_j > 1$ for some $j \in [d]$.  Since $P_\sigma$ is bounded, we
    can find $i \in I_\sigma$ with $(a_i)_j \neq 0$. Since $\sigma \a_i \in
    \Z^d_{\ge0}$, we compute $1 \ge \inner{\a_i,\vb} \ge (a_i)_j v_j > 1$,
    which is a contradiction. Hence all facets of $P_\sigma$ have facet width
    one and, by Theorem~\ref{thm:compressed}, $P_\sigma$ is compressed.
\end{proof}

\begin{theorem}
\label{thm:LAB-triangulation}
    If $P$ is a reflexive and locally anti-blocking polytope, then $P$ has a
    regular and unimodular triangulation. In particular, $h^*(P)$ is unimodal.
\end{theorem}
\begin{proof}
    By definition of locally anti-blocking polytopes, $P$ is subdivided into
    the polytopes $P_\sigma$ for $\sigma \in \PM^d$. As is easily seen, this
    is a regular subdivision with respect to the height function $\omega(\vb)
    := \|\vb\|_1$. 
    
    Let $U := P \cap \Z^d$.  By Lemma~\ref{lem:LAB-reflexive}, the polytopes
    $P_\sigma$ are compressed and thus, using Theorem~\ref{thm:compressed}, $U
    \cap \R^d_\sigma$ is the vertex set of $P_\sigma$. Fix a order on $U$.
    Since $P_\sigma$ is compressed, the pulling triangulation $\Triang_\sigma$
    of $P_\sigma$ is unimodular and we only need to argue that $\Triang =
    \bigcup_\sigma \Triang_\sigma$ is a regular triangulation of $P$.  Let
    $\sigma, \sigma' \in \PM^d$ and let $S \in \Triang_\sigma$ and $S' \in
    \Triang_{\sigma'}$ be two simplices.  Then $S \cap S'$ is contained in $F
    = P_\sigma \cap P_{\sigma'}$, which is a face of $P_\sigma$ as well as of
    $P_{\sigma'}$. By the construction of a pulling triangulation, we see that
    both $S \cap F$ and $S' \cap F$ are simplices of the pulling triangulation
    of $F$ and hence $S \cap S' = (S \cap F) \cap (S' \cap F)$ is a face of
    both $S$ and $S'$.  The same argument as in the proof of
    Theorem~\ref{thm:ExtendTriang}, then shows that $\Triang$ is a regular
    triangulation. The unimodality of $h^*(P)$ now follows from
    Theorem~\ref{thm:BrunsRoemer}.
\end{proof}

\begin{remark} \label{rmk:generalIDP}
The techniques of this section can be extended to the following class of
polytopes. We say that a polytope $P \subset \R^d$ has the
\Def{orthant-lattice property} (OLP) if the restriction $P_\sigma = P \cap
\R^d_\sigma$ is a (possibly empty) lattice polytope for all $\sigma \in
\PM^d$. If $P$ is reflexive, then $P_\sigma$ is full-dimensional for every
$\sigma$. Now, if every $P_\sigma$ has a unimodular cover, then so does $P$
and hence is IDP.  Let $P_\sigma = \{ \x \in \R^d_\sigma : A^\sigma \x \le
b^\sigma\}$. Then some conditions that imply the existence of a unimodular
cover include:
\begin{enumerate}
		\item $P_\sigma$ is compressed;
		\item ${A}^\sigma$ is a totally unimodular matrix;
		\item ${A}^\sigma$ consists of rows which are $B_d$ roots;
		\item $P_\sigma$ is the product of unimodular simplices;
		\item There exists a projection $\pi:\R^d\to \R^{d-1}$ such that
            $\pi(P_\sigma)$ has a regular, unimodular triangulation
            $\Triang$ such that the pullback subdivision $\pi^\ast
            (\Triang)$ is lattice.
	\end{enumerate}
We refer to~\cite{HPPS} for background and details.

An example of such a  polytope is
\[
	P=\conv \begin{bmatrix}
	1 & 0 & 0 & 1 & 0 & 1 & -1 & 0 & 0 \\
	0 & 1 & 0 & 1 & 1 & 1 & 0 & -1 & 0\\
	0 & 0 & 1 & 0 & 1 & 1 & 0 & 0 & -1\\
	\end{bmatrix}\subset \R^3.
\]
This is a reflexive OLP polytope. The restriction to $\R^3_+$ is 
	\[
	P_+=\conv \begin{bmatrix}
	1 & 0 & 0 & 1 & 0 & 1 & 0\\
	0 & 1 & 0 & 1 & 1 & 1 & 0\\
	0 & 0 & 1 & 0 & 1 & 1 & 0
	\end{bmatrix}\subset \R^3 \, ,
	\]
which is not an anti-blocking polytope.
\end{remark}

The \Def{Mahler conjecture} in convex geometry states that every
centrally-symmetric convex body $K \subset \R^d$ satisfies
\[
    \vol(K) \cdot \vol(K^*) \ \ge \ 
    \vol(C_d) \cdot \vol(C_d^*) \, ,
\]
where $C_d = [-1,1]^d$ is the $d$-cube. The Mahler conjecture has been
verified only in small dimensions and for special classes of convex bodies. In
particular, Saint-Raymond~\cite{Saint-Raymond} proved the following beautiful
inequality, where $A(P)$ refers to the anti-blocking dual of $P$, see Section~\ref{sec:graphs}. The characterization of the equality case is independently due
to Meyer~\cite{Meyer-unconditional}  and Reisner~\cite{Reisner-unconditional}.

\begin{theorem}[Saint-Raymond]\label{thm:raymond}
    Let $P \subset \R^d_+$ be an anti-blocking polytope. Then
    \[
        \vol(P) \cdot \vol(A(P)) \ \ge \ \frac{1}{d!}
    \]
     with equality if and only if $P$ or $A(P)$ is the cube $[0,1]^d$.
\end{theorem}

This inequality directly implies the Mahler conjecture for unconditional
reflexive polytopes, that we record for the normalized volume.

\begin{corollary} \label{cor:mahler}
    Let $P \subset \R^d$ be an unconditional reflexive polytope. Then
    \[
        \Vol(P) \cdot \Vol(P^*) \ \ge \ 4^d d! 
    \]
    with equality if and only if $P$ or $P^*$ is the cube $[-1,1]^d$.
\end{corollary}

\section{Unconditional reflexive polytopes and perfect
graphs}\label{sec:graphs}

For $A \subseteq [d]$, let $\1_A \in \{0,1\}^d$ be its characteristic vector.
If $\Gamma \subseteq 2^{[d]}$ is a simplicial complex, that is, a nonempty set
system closed under taking subsets, then 
\[
    P \ = \ \conv( \1_\sigma : \sigma \in \Gamma)
\]
is an anti-blocking $0/1$-polytope and every anti-blocking polytope with
vertices in $\{0,1\}^d$ arises that way (cf.~\cite[Thm.~2.3]{Fulkerson2}). 

A prominent class of anti-blocking $0/1$-polytopes arises from graphs. Given a graph $G=([d],E)$ with $E \subseteq \binom{[d]}{2}$, we say that
$S\subseteq [d]$ is a \Def{stable set} (or \Def{independent set}) of $G$ if $uv
\not\in E$ for any $u,v \in S$.  The \Def{stable set polytope} of $G$ is
\[
    P_{G} \ \coloneqq \ \conv\{\1_S \ : \ S \subseteq [d] \text{ stable set of
    $G$} \} \,
    .
\] 
Since the stable sets of a graph form a simplicial complex, $P_G$ is an
anti-blocking polytope.  Stable set polytopes played an important role in the
proof of the \emph{weak perfect graph conjecture}~\cite{Lovasz}. A
\Def{clique} is a set $C\subseteq [d]$ such that every two vertices in $C$ are
joined by an edge.  The \Def{clique number} $\omega(G)$ is the largest size of
a clique in $G$. A graph is \Def{perfect} if $\omega(H)= \chi (H)$ for all
induced subgraphs $H\subseteq G$, where $\chi(H)$ is the chromatic
number of $H$.

Lov\'asz~\cite{Lovasz} gave the following geometric characterization of
perfect graphs; see also~\cite[Thm.~9.2.4]{GLS}.  For a set $C \subseteq [d]$
and $\x \in \R^d$, we write $\x(C) := \inner{\1_C,\x} = \sum_{i \in C} x_i$.
\begin{theorem} \label{thm:Lovasz}
    A graph $G=([d],E)$ is perfect if and only if 
	\[
        P_G \ = \ \{ \x\in\R^d_+  :  \ \x(C) \leq 1
        \text{ for all cliques }  C\subseteq [d] \}.
	\]
\end{theorem} 	

For an anti-blocking polytope $P \subset \R^d_+$ define the \Def{anti-blocking
dual}
\[
    A(P) \ \coloneqq \  \{ \y \in \R^d_+ : \inner{\y,\x} \le 1 \text{ for all } \x
    \in P\} \, .
\]
The polar dual $(\U P)^*$ is again unconditional and it follows that 
\[
    (\U P)^* \ = \ \U A(P) \, .
\]

\begin{theorem}[{\cite[Thm.~9.4]{schrijver1986}}] \label{thm:AB-polar}
    Let $P \subset \R^d_+$ be a full-dimensional anti-blocking polytope with
    \begin{align*}
        P \ = \ \convDown{\c_1,\dots,\c_r} &\ = \ \{ \x \in \R^n_+ :
        \inner{\d_i,\x} \le 1 \text{ for all } i=1,\dots,s\} \\
    \intertext{for some $\c_1,\dots,\c_r,\d_1,\dots,\d_s \in \R^d_+$.  Then}
        A(P) \ = \ \convDown{\d_1,\dots,\d_s} &\ = \ \{ \x \in \R^d_+ : 
         \inner{\c_i,\x} \le 1 \text{ for all } i=1,\dots,r\}.
    \end{align*}
    In particular, $A(A(P)) = P$.
\end{theorem}

Theorems~\ref{thm:Lovasz} and~\ref{thm:AB-polar} then imply for a perfect
graph $G$ that
\begin{equation}\label{eqn:complement}
    A(P_G) \ = \ \convDown{\1_C : C \text{ clique in } G } \ = \
    P_{\overline{G}} \, ,
\end{equation}
where $\overline{G} = ([d], \binom{[d]}{2} \setminus E)$ is the
\Def{complement graph}.

\begin{corollary}[Weak perfect graph theorem]
    A graph $G$ is perfect if and only if $\overline G$ is perfect.
\end{corollary}

We note that in particular if $G$ is perfect, then $P_G$ is compressed.

\begin{proposition}[{\cite[Prop.~3.10]{ChappellFriedlSanyal}}]
    \label{prop:AB-compressed}
    Let $P \subset \R^d_+$ be an anti-blocking polytope. Then $P$ is
    compressed if and only if $P = P_G$ for some perfect
    graph $G$.
\end{proposition}

Let us remark that Theorem~\ref{thm:Lovasz} also allows us to characterize the
Gorenstein stable set polytopes. For comparability graphs of posets (see
Section~\ref{sec:ChainAndGroebner}) this was noted by
Hibi~\cite{Hibi-LatticesAndStraighteningLaws}. A graph $G$ is called
\Def{well-covered} if every inclusion-maximal stable set has the same size. It
is called \Def{co-well-covered} if $\OV G$ is well-covered.

\begin{proposition}[{\cite[Theorem 2.1b]{OhsugiHibi-Conjecture}}] \label{prop:StableGorenstein}
    Let $P_G$ be the stable set polytope of a perfect graph $G= ([d], E)$.
    Then $P_G$ is Gorenstein if and only if $G$ is co-well-covered.
\end{proposition}
\begin{proof}
    By definition, $P_G$ is Gorenstein if there are $c \in \Z_{>0}$ and $\q \in
    \Z^d$ such that $\q + c\, P_G$ is reflexive. Using
    Theorem~\ref{thm:Lovasz}, we see that $\q + c\, P_G$ is given by all points
    $\x \in \R^d$ such that
    \[
        -x_i \ \le -\ q_i \quad \text{ for } i = 1,\dots,d 
        \qquad \text{ and } \qquad
        \x(C) \ \le c + \q(C) \quad \text{ for all maximal cliques } C
        \subseteq [d] 
    \]
    These inequalities are facet-defining, as can be easily seen. Using
    the representation given in Definition~\ref{def:reflexive}, we note that
    $\q + c\, P_G$ is reflexive if and only if $q_i = -1$ for all
    $i=1,\dots,d$ and $c + \q(C) = c - |C| = 1$ for all maximal cliques $C$.
    This happens if and only if all maximal cliques have the same size.
\end{proof}

Combining Lemma~\ref{lem:LAB-reflexive} with
Proposition~\ref{prop:AB-compressed} yields the following characterization of
reflexive locally anti-blocking polytopes.

\begin{theorem} \label{thm:unconditional}
    Let $P\subset \R^d$ be a locally anti-blocking lattice polytope with $0$
    in its interior. Then $P$ is reflexive if and only if for every $\sigma
    \in \PM^d$ there is a perfect graph $G_\sigma$ such that $P_\sigma =
    \sigma P_{G_\sigma}$.

    In particular, $P$ is an unconditional reflexive polytope if and only
    if $P = \U P_G$ for some perfect graph $G$.
\end{theorem}

The following corollary to Theorem~\ref{thm:unconditional} was noted in
\cite[Thm.~3.4]{ChappellFriedlSanyal}. The second part also appears
in~\cite[Example 2.3]{OH18}.

\begin{corollary}
    If $G_1,G_2$ are perfect graphs on the vertex set $[d]$, then $P_{G_1} +
    (-P_{G_2})$ and $P_{G_1} \vee (-P_{G_2})$ are reflexive polytopes.
\end{corollary}

For $G_1 = G_2 = K_d$ the complete graph on $d$ vertices, the polytope
$P_{G_1} + (-P_{G_2})$ is the Legendre polytope studied by
Hetyei et al.~\cite{hetyei,EHR}.

Using \textsc{Normaliz}~\cite{Normaliz} and the Kreuzer--Skarke database for
reflexive polytopes~\cite{kreuzerskarke,kreuzerskarke4dim}, we were able to
verify that $72$ of the $3$-dimensional reflexive polytopes and at least $407$
of the $4$-dimensional reflexive polytopes with at most $12$ vertices are
locally anti-blocking.  Unfortunately, our computational resources were too
limited to test \emph{most} of the $4$-dimensional polytopes. However, there
are only $11$ $4$-dimensional unconditional reflexive polytopes (by virtue of
Theorem~\ref{thm:uncond_refl_different}).

If $G, G'$ are perfect graphs, then $G \uplus G'$ as well as its bipartite sum
$G \bowtie G' = \overline{ \overline{G} \uplus \overline{G'}}$ are perfect.
On the level of unconditional polytopes we note that
\[
    \U P_{G \uplus G'} \ = \ \U P_G \times \U P_{G'} \quad \text{ and } \quad \U
    P_{G \bowtie G'} \ = \ \U P_G \oplus \U P_{G'} \, ,
\]
where $\oplus$ is the \Def{direct sum} (or \Def{free sum}) of
polytopes~\cite{Henk}.  These observations give us the class of Hanner
polytopes which are important in relation to the $3^d$-conjecture;
see~\cite{SWZ}. A centrally symmetric polytope $H \subset \R^d$ is called a
\Def{Hanner polytope} if and only if $H = [-1,1]$ or $H$ is of the form $H_1
\times H_2$ or $H_1 \oplus H_2 = (H_1^* \times H_2^*)^*$ for lower dimensional
Hanner polytopes $H_1,H_2$.  Thus, every Hanner polytope is of the form $\U
P_G$ for some perfect graph $G$. Hanner polytopes were obtained from split
graphs in~\cite{FHSZ} using a different geometric construction.

Let us briefly note that Theorem~\ref{thm:unconditional} also yields bounds on
the entries of the $h^*$-vector. Recall that $h^*_i(C_d)$ for the cube $C_d =
[-1,+1]^d$ is given by the \Def{type-B Eulerian number} $B(d,i) =
\sum_{j=1}^{i} (-1)^{i-j}\binom{d}{i-j}(2j-1)^{d-1}$, that counts signed
permutations with $i$ descents (see also Section~\ref{sec:BBn}). Its polar
$C_d^*$ is the \Def{crosspolytope} with $h^*_i(C_d^*) = \binom{d}{i}$ for
$i=0,\dots,d$.

\begin{corollary}\label{cor:h-bounds}
    Let $P \subset \R^d$ be an unconditional reflexive polytope. Then 
    \[
        \binom{d}{i} \ \le \ h^*_i(P) \ \le \ B(d,i) \, .
    \]
\end{corollary}
\begin{proof}
    It follows from Theorem~\ref{thm:unconditional} that every reflexive and
    unconditional $P$ satisfies $C_d^* \subseteq P \subseteq C_d$, where $C_d
    = [-1,1]^d$. By Theorem~\ref{thm:h-mono}, the entries of the
    $h^*$-vector are monotone with respect to inclusion.
\end{proof}

We close the section by showing that distinct perfect graphs yield distinct
unconditional reflexive polytopes.

\begin{theorem}\label{thm:uncond_refl_different}
    Let $G,H$ be perfect graphs on vertices $[d]$. Then $\U P_G$ is unimodularly
    equivalent to $\U P_H$ if and only if $G \cong H$.
\end{theorem}
\begin{proof}
    Assume that $T(\U P_G) = \U P_H$ for some $T(\x) = W\x + \t$ with $\t
    \in \Z^d$ and $W \in \operatorname{SL}_d(\Z)$. Since the origin is the
    only interior lattice point of both polytopes, we infer that $\t = 0$.
    Let $W = (\w_1,\dots,\w_d)$. Thus, $\z \in \Z^d$ is a lattice point in $
    \U P_H$ if and only if there is a stable set $S$ and $\sigma \in \PM^S$
    such that 
    \begin{equation}\label{eqn:sW}
        \z \ = \ \sum_{i \in S} \sigma_i \w_i.
    \end{equation}

    On the one hand, this implies that $\w_i$ and $\w_j$ have disjoint
    supports whenever $i,j \in S$ and $i \neq j$. Indeed, if the supports of
    $\w_i$ and $\w_j$ are not disjoint, then $\sigma_i \w_i + \sigma_j \w_j$
    has a coordinate $>1$ for some choice of $\sigma_i,\sigma_j \in \PM$,
    which contradicts the fact that $\U P_H \subseteq [-1,1]^d$.

    On the other hand, for any $h \in [d]$, the point $\e_h = \1_{\{h\}}$ is
    contained in $\U P_H$. Hence, there is a stable set $S$ and $\sigma \in
    \PM^S$ such that~\eqref{eqn:sW} holds for $\z = \e_h$. Since the
    supports of the vectors indexed by $S$ are disjoint, this means that $S
    = \{ i \}$ and $\e_h = \sigma_i \w_i$. We conclude that $W$ is a signed
    permutation matrix and $G \cong H$.
\end{proof}
 
We can conclude that the number of unconditional reflexive polytopes in $\R^d$
up to unimodular equivalence is precisely the number of unlabeled perfect
graphs on $d$ vertices. This number has been computed up to $d=13$ (see
\cite[Sec.5]{Hougardy-PerfectGraphs} and
\href{https://oeis.org/A052431}{A052431} of \cite{oeis}). We show the
sequence in Table~\ref{table:PerfectGraphs}.
  \begin{table}
 	 \begin{tabular}{|c|c|c|c|c|c|c|c|c|c|c|c|c|c|}
 	\hline 
 	n &  3 & 4 & 5 & 6 &7 &8 &9 &10 &11 & 12
     & 13 \\ \hline
     p(n) &  4 & 11 & 33 & \smaller[1]148 & \smaller[1]906 & \smaller[2] 8887 & \smaller[2]136756 &
     \smaller[2] 3269264 & \smaller[2] 115811998& \smaller[2] 5855499195 &
     \smaller[2] 410580177259 \\ \hline 
 	\end{tabular}
     \caption{Number $p(n)$ of unlabeled perfect graphs; OEIS sequence \href{https://oeis.org/A052431}{A052431}.}
 \label{table:PerfectGraphs}
 \end{table}

\section{The type-$B$ Birkhoff polytope} \label{sec:BBn}
The \Def{Birkhoff polytope} $\mathcal{B}(n)$ is defined as the convex hull of
all $n\times n$ permutation matrices. Equivalently, $\mathcal{B}(n)$ is the
set of all \Def{doubly stochastic matrices}, that is, nonnegative matrices $M$
with row and column sums equal to $1$, by work of Birkhoff \cite{Birkhoff}
and, independently, von Neumann \cite{vonNeumann}.  This polytope has been
studied quite extensively and is known to have many properties of interest
(see, e.g.,
\cite{athanasiadismagic,Bapat,beckpixton,CanfieldMcKay,Davis-stochastic,DeLoeraLiuYoshida,Paffenholz-Birkhoff}).
Of particular interest to our purposes, it is known to be Gorenstein, to be
compressed \cite{StanleyNonnegativity}, and to be $h^\ast$-unimodal
\cite{athanasiadismagic}.  In this section, we will introduce a type-$B$
analogue of this polytope corresponding to signed permutation matrices and
verify many similar properties already known for $\mathcal{B}(n)$.

The \Def{hyperoctahedral group} is defined to by $B_n \coloneqq \left(\Z/2\Z \right) \wr \mathfrak{S}_n$, which is the Coxeter group of type-$B$ (or type-$C$).
Elements of this group can be thought of as permutations from $\mathfrak{S}_n$ expressed in one-line notation $\sigma=\sigma_1\sigma_2\cdots\sigma_n$, where we also associate a sign $\sgn(\sigma_i)$ to each $\sigma_i$. 
To each signed permutation $\sigma\in B_n$, we associate a matrix $M_\sigma$ defined as $\left(M_\sigma\right)_{i,\sigma_i}=\sgn(\sigma_i)$ and $\left(M_\sigma\right)_{i,j}=0$ otherwise. If every entry of $\sigma$ is positive, then $M_{\sigma}$ is simply a permutation matrix.
This leads to the following definition:

\begin{definition}
The \Def{type-$\boldsymbol B$ Birkhoff polytope} (or \Def{signed Birkhoff polytope}) is
	\[
	\BB(n)\coloneqq\conv\left\lbrace M_\sigma \, : \, \sigma\in B_n \right\}\subset \R^{n\times n}.
	\]
That is, $\BB(n)$ is the convex hull of all $n\times n$ signed permutation matrices. 
\end{definition}

This polytope was previously studied in \cite{ConvHullCoxeter}, though the
emphasis was not on Ehrhart-theoretic questions. Since all points in the
definition of $\BB(n)$ lie on a sphere, it follows that they are all vertices.

\begin{proposition}
    For every $\sigma\in B_n$, $M_\sigma$ is a vertex of $\BB(n)$.
\end{proposition}

It is clear that $\BB(n)$ is an unconditional lattice polytope in $\R^{d
\times d}$ and we study it by restriction to the positive orthant.

\begin{definition}
    For $n \ge 1$, we define the \Def{positive type-$B$ Birkhoff polytope},
    $\pos(n)$, to be the polytope
    \[
        \pos(n) \coloneqq \BB(n) \cap \R^{n \times n}_+.
    \]
\end{definition}

A simple way to view this as an anti-blocking polytope is via matching
polytopes.  Given a graph $G=([d],E)$, a \Def{matching} is a set $M \subseteq
E$ such that $e \cap e' = \emptyset$ for any two distinct $e,e' \in M$. The
corresponding \Def{matching polytope} is
\[
    \Mat(G) \ \coloneqq \ \conv\{\1_M : M\subseteq E \text{ matching}\}
    \subset \R^{E} \, .
\]
If $G$ is a bipartite graph, then the matching polytope is easy to describe.
For $v \in [d]$ let $\delta(v) \subseteq E$ denote the edges incident to $v$. 

\begin{theorem}[{\cite[Sec. 8.11]{schrijver1986}}]\label{thm:match}
    For bipartite graphs $G$ the matching polytope is given by
    \[
        \Mat(G) \ = \ \{ \x \in \R^E_+ : \x(\delta(v)) \le 1 \text{ for all } v
        \in [d] \} \, .
    \]
\end{theorem}

As a simple consequence, we get 
\begin{corollary}
    $\pos(n)$ is the matching polytope of the complete bipartite graph
    $K_{n,n}$ on $2n$ vertices. 
\end{corollary}
\begin{proof}
    We can identify the edges of $K_{n,n}$ with $[n] \times [n]$.  Every
    matrix $M \in \pos(n) \cap \{0,1\}^{n \times n}$ is a partial permutation
    matrix and therefore contains at most one $1$ in every row and column.  It
    follows that the set $\{ (i,j) : M_{ij} = 1 \}$ is a matching of $K_{n,n}$
    and every matching arises that way. Since $\pos(n)$ and 
    $\Mat(K_{n,n})$ are both $0/1$-polytopes, this proves the claim.
\end{proof}

It follows from the description given in Theorem~\ref{thm:match} and the
definition of compressed polytopes that matching polytopes of bipartite graphs
are compressed.  Hence, by Proposition~\ref{prop:AB-compressed} $\Mat(G)$ is
the stable set polytope of a perfect graph. The graph in question is the
\Def{line graph} $L(G)$ on the vertex set $E$ and edge $ee'$ whenever $e\cap
e' \neq \emptyset$. It is clear that $M$ is a matching in $G$ if and only if
$M$ is a stable set in $L(G)$. If $L(G)$ is perfect, then $G$ is called a
\Def{line perfect graph}.  From Lov\'asz' Theorem~\ref{thm:Lovasz} one can then
infer $\Mat(G) = P_{L(G)}$ and hence bipartite graphs are line perfect;
cf.~\cite[Thm.  2]{Maffray-LineGraphs}.

The polytope $\pos(n)$ is the stable set polytope of $L(K_{n,n})=K_n\square
K_n$, the Cartesian product of complete graphs, which is the graph of
legal moves of a rook on an $n$-by-$n$ chessboard 
and thus called a \Def{rook graph}.

Since all vertices in $K_{n,n}$ have the same degree, it follows that all
maximal cliques in $K_n \square K_n$ have size $n$ and from
Proposition~\ref{prop:StableGorenstein} we conclude the following.

\begin{corollary} \label{cor:PosGorenstein}
    The polytope $\pos(n)$ is Gorenstein.
\end{corollary}

For two matrices $A,B \in \R^{d \times d}$ we denote by $\inner{A,B} =
\mathrm{tr}(A^tB)$ the Frobenius inner product. Also, for vectors $\u,\vb \in
\R^d$ let us write $ \u \otimes \vb \in \R^{d \times d}$ for the matrix with
$(\u \otimes \vb)_{ij} = u_i v_j$.

\begin{corollary}
    The polytope $\BB(n)$ is an unconditional reflexive polytope. Its
    facet-defining inequalities are given \rem{by $A_{ij} \ge 0$ and}  by
    \[
        \inner{A, \sigma \otimes \e_i } \ \le \ 1 \quad \text{ and } \quad
        \inner{A, \e_i \otimes \sigma } \ \le \ 1
    \]
    for all $i=1,\dots,n$ and $\sigma \in \PM^n$.
\end{corollary}
The inequality description of this polytope was previously obtained
in~\cite{ConvHullCoxeter} using the notion of \emph{Birkhoff tensors}.

\begin{proof}
    We deduce that $\BB(n)$ is reflexive by appealing to Theorem
    \ref{thm:unconditional}, using the fact that $\BB(n)$ is unconditional and
    $\pos(n) = \BB(n) \cap \R^{n \times n}_+$ is the stable set polytope of a
    perfect graph as discussed above. We obtain the inequality description by
    applying Proposition~\ref{prop:AB-inequalities} and
    Theorem~\ref{thm:match}.
\end{proof}

The dual $\Cn \coloneqq \BB(n)^*$ is the unconditional reflexive polytope
associated with the graph $\overline{K_n\square K_n}$. The corresponding
anti-blocking polytope $\Cposn=P_{\overline{K_n\square K_n}}$ also has the
nice property that all cliques have the same size $n$ and hence
Proposition~\ref{prop:StableGorenstein} applies.

\begin{corollary}
    The polytope $\Cposn$ is Gorenstein.
\end{corollary}

By Theorems~\ref{thm:LAB-triangulation}, Theorem~\ref{thm:unconditional}, and
Proposition \ref{prop:AB-compressed}, we have the following unimodality
results.
\begin{corollary}
    For any $n\in\Zpos$, we have that $h^\ast(\BB(n))$, $h^\ast(\pos(n))$,
    $h^\ast(\Cn)$, and $h^\ast(\Cposn)$ are unimodal.
\end{corollary}

Let us conclude this section with some enumerative data. The polytope $\BB(n)$
has $2^n n!$ vertices and $n 2^{n+1}$ facets.  In contrast, the vertices of
$\pos(n)$ are in bijection to partial permutations of $[n]$. Hence $\pos(n)$
has  $n! \sum_{i=0}^n \frac{1}{i!}$ many vertices but only $n^2 + 2n$ facets.
The polytope $\Cposn$ has $n 2^{n+1} - (n+1)^2$ many vertices and $n^2 + n!$
facets. We used \textsc{Normaliz}~\cite{Normaliz} to compute the normalized
volume and $h^*$-vectors of these polytopes; see Tables~\ref{table:pos},
\ref{table:BBn}, \ref{table:posdual}, and \ref{table:dual}.  Given the
dimension and volumes of these polytopes, our computational resources were
quite quickly exhausted. Note that $\BB(3)$ and $\mathcal{C}(3)$ have
precisely the same Ehrhart data and normalized volume and in fact it is
straightforward to verify that  $\BB(3)$ and $\mathcal{C}(3)$ are unimodularly
equivalent.

Using Theorem~\ref{thm:raymond} and Corollary~\ref{cor:mahler}, we get a lower
bound on the volume of $\pos(5)$ and $\BB(5)$, respectively. We get that 
\[
\begin{array}{l@{\quad > \quad }r}
\Vol(\BB_+(5))  &  30.637.007.047.800 \\
\Vol(\BB(5))    &  1.028.007.369.668.940.603.880\\
\end{array}
\]
are bounds on the number of simplices in an unimodular triangulation.

\begin{center} 
\begin{table}
	\begin{tabular}{|c||c|c|}
	\hline
	$n$ & $\operatorname{Vol}(\pos(n))$ & $h^\ast(\pos(n))$\\ \hline
	$1$ & $1$ & $1$\\ \hline
	$2$ & $4$ & $(1,2,1)$\\ \hline
	$3$ & $642$ & $(1,24,156,280,156,24,1)$\\ \hline 
	$4$ & $12065248$ & $(1,192,9534,151856, 975793,2860752,4069012, $\\ 
	       &                       & $2860752,975793,151856,9534,192,1)$\\\hline
	\end{tabular}
\caption{$\pos(n)$.}
\label{table:pos}	
\end{table}
\end{center}

\begin{center}
\begin{table}
	\begin{tabular}{|c||c|c|}
	\hline 
	$n$ & $\operatorname{Vol}(\BB(n))$ & $h^\ast(\BB(n))$\\ \hline
	$1$ & $2$ & $(1,1)$\\ \hline 
	$2$ & $64$ & $(1,12,38,12,1)$\\ \hline
	$3$ & $328704$ & $(1,129,4482,40844,118950,118950,40844,4482,129,1)$\\ \hline
	$4$ & $790708092928$ & ? \\ \hline
	\end{tabular}
\caption{$\BB(n)$.}
\label{table:BBn}
\end{table}
\end{center}

\begin{center} 
\begin{table}
	\begin{tabular}{|c||c|c|}
	\hline
	$n$ & $\operatorname{Vol}(\Cposn)$ & $h^\ast(\Cposn)$\\ \hline
	$1$ & $1$ & $1$\\ \hline
	$2$ & $6$ & $(1,4,1)$\\ \hline
	$3$ & $642$ & $(1, 24, 156, 280, 156, 24, 1)$\\ \hline 
	$4$ & $2389248$ & $(1, 88, 2656, 34568, 201215, 562112, 787968  $\\ 
	       &                       & $562112, 201215, 34568, 2656, 88, 1 )$\\ \hline
	 $5$ & $506289991680$ & $?$\\ \hline 
	\end{tabular}
\caption{$\Cposn$.}
\label{table:posdual}	
\end{table}
\end{center}

\begin{center}
\begin{table}
	\begin{tabular}{|c||c|c|}
	\hline 
	$n$ & $\operatorname{Vol}(\Cn)$ & $h^\ast(\Cn)$\\ \hline
	$1$ & $2$ & $(1,1)$\\ \hline 
	$2$ & $96$ & $(1, 20, 54, 20, 1 )$\\ \hline
	$3$ & $328704$ & $(1, 129, 4428, 40844, 118950, 118950, 40844, 4428, 129, 1 )$\\ \hline
	$4$ & $156581756928$ & $(1, 592, 110136, 8093168, 222332060, 2558902352, $ \\
	       &                              & $13699272072, 36553260912, 50497814342, 36553260912, $\\
	       &                              &$13699272072 2558902352, 222332060, 8093168, 110136, 592, 1)$ \\ \hline
	$5$ & $16988273098107125760$ & $?$\\ \hline 
	\end{tabular}
\caption{$\Cn$.}
\label{table:dual}
\end{table}
\end{center}

\section{CIS graphs and compressed Gale-dual pairs of
polytopes}\label{sec:gale}
The notion of Gale-dual pairs was introduced in \cite{GardnerPolytope}. Given
two polytopes $P,Q\subset \R^d$, we say that these polytopes form a
\Def{Gale-dual pair} if
\begin{align*}
    P \ &= \ \{\x\in\R^d_+  : \ \inner{\x,\y} =1 \text{ for }
    \y \in V(Q)\}  \text{ and } \\
    Q \ &= \ \{\x\in\R^d_+  : \ \inner{\x,\y} =1 \text{ for }
    \y \in V(P)\} \, .
\end{align*}

The prime example of a Gale-dual pair of polytopes is the Birkhoff polytope
$\mathcal{B}_n$, the convex hull of permutation matrices $M_\tau$, and the
\Def{Gardner polytope} $\mathcal{G}_n$, which is the polytope of all
nonnegative matrices $A \in \R^{n \times n}_+$ such that $\inner{M_\tau,A} =
1$ for all permutation matrices $M_\tau$. Both polytopes are compressed,
Gorenstein lattice polytopes of codegree $n$. The question raised
in~\cite{GardnerPolytope} was if there were other Gale-dual pairs with (a subset
of) these properties. In this section we briefly outline a construction for
compressed Gale-dual pairs of polytopes.

Following~\cite{ABG}, we call $G = ([d],E)$ a \Def{CIS graph} if $C \cap S
\neq \emptyset$ for every inclusion-maximal clique $C$ and inclusion-maximal
stable set $S$.  For brevity, we refer to those as maximal cliques and stable
sets, respectively.  For example, if $B$ is a bipartite graph with \Def{perfect
matching}, i.e., a matching covering all vertices, then the line graph $L(G)$
is CIS. Another class of examples is given by a theorem of
Grillet~\cite{Grillet}. Let $\Pi = ([d],\preceq)$ be a partially ordered set.
The \Def{comparability graph} of $\Pi$ is the simple graph $G_\prec = ([d],E)$
with $ij \in E$ if $i \prec j$ or $j \prec i$.  Comparability graphs are known
to be perfect~\cite{Hougardy-PerfectGraphs}.  The \Def{bull graph} is the
graph vertices $a,b,c,d,e$ and edges $ab, bc, cd, de, bd$.

\begin{theorem}[\cite{Grillet}]\label{thm:Grillet} 
    Let $(\Pi,\preceq)$ be a poset with comparability graph $G$. Then $G$ is
    CIS if every induced $4$-path is contained in an induced bull graph.
\end{theorem}

The wording in graph-theoretic terms is due to Berge; see~\cite{Zang} for
extensions. 

\begin{proposition}\label{prop:CIS}
    Let $G$ be a perfect CIS graph. Then 
    \begin{align*}
        P \ &= \ \conv( \1_S : S \text{ maximal stable set of } G  )\\ 
        Q \ &= \ \conv( \1_C : C \text{ maximal clique of } G  )
    \end{align*}
    is a Gale-dual pair of compressed polytopes.
\end{proposition}
\begin{proof}
    For a stable set $S$ and a clique $C$, we have that $S \cap C \neq \emptyset$
    if and only if $\inner{\1_C,\1_S} = 1$. Let $P_G$ be the stable set
    polytope of $G$. It follows from Theorem~\ref{thm:Lovasz} that the
    vertices of 
    \[
        P' \ := \ P_G \cap \bigcap_{C \text{ maximal clique}} \{ \x \in \R_+^n
        : \x(C) = \inner{\1_C,\x} = 1 \} 
    \]
    are of the form $\1_S$ for stable sets meeting every maximal clique
    non-trivially. Note that a stable set $S$ of the CIS graph $G$ is maximal
    if $S \cap C \neq \emptyset$ for every maximal clique $C$. Hence $P = P'$ and $P$ is a face of $P_G$. Since faces
    of compressed polytopes are compressed, it follows that $P$ is compressed.

    The complement graph $\OV G$ is also a perfect CIS graph and the same
    argument applied to $Q$ completes the proof.
\end{proof}

Note that both of the examples above are perfect CIS graphs. This shows
that compressed (lattice) Gale-dual pairs are not rare. Recall that a
graph $G$ is well-covered if every maximal stable set has the same size and
$G$ is co-well-covered if $\OV G$ is well-covered.
Theorem~\ref{thm:Grillet} and its generalization in~\cite{Zang} allow for the
construction of perfect CIS graphs which are well-covered and co-well-covered
(for example, by taking ordinal sums of antichains). Moreover, the recent
paper~\cite{DHMV} gives classes of examples of well-covered and
co-well-covered CIS graphs. This is a potential source of compressed
Gorenstein Gale-dual pairs but we were not able to identify the perfect
graphs in these families.

Theorem~\ref{thm:unconditional} implies that if $(F,G)$ is a Gale-dual pair of
Proposition~\ref{prop:CIS}, then there is a (unconditional) reflexive polytope
such that $F \subset P$ and $G \subset P^*$ are dual faces.

\begin{question}
    Is it true that every Gale-dual pair $(F,G)$ appears as dual faces of some
    reflexive polytope $P$?
\end{question}

\section{Chain Polytopes and Gr\"obner Bases}
\label{sec:ChainAndGroebner}

Given a lattice polytope $P \subset \R^d$, the existence of regular
triangulations, particularly those which are unimodular and flag, has direct
applications to the associated toric ideal of $P$.  In this section, we will
discuss how certain Gr\"obner bases of the toric ideal of an anti-blocking
polytope can be extended to Gr\"obner bases of the associated unconditional
polytope.  In particular, we provide an explicit description of Gr\"obner
bases for unconditional polytopes arising from the special class of
anti-blocking polytopes called chain polytopes. We refer the reader to the
wonderful books~\cite{CLO} and~\cite{Sturmfel-GrobnerPolytopes} for background
on Gr\"obner bases and toric ideals.

Let $Z \coloneqq P \cap \Z^d$.  The \Def{toric ideal associated to $P$} is the
ideal $I_P \subset \C[ x_{\p} : \p \in Z]$ with generators
\[
    x_{\r_1}x_{\r_2} \cdots x_{\r_k} - x_{\s_1}x_{\s_2} \cdots x_{\s_k} \, ,
\]
where $\r_1,\dots,\r_k,\s_1,\dots,\s_k \in Z$ are lattice points such that
$\r_1 + \cdots + \r_k = \s_1 + \cdots + \s_k$. If we denote the two
\emph{multi}sets of points by $R$ and $S$, we simply write $x^R - x^S$.  A
celebrated result of Sturmfels~\cite[Thm.~8.3]{Sturmfel-GrobnerPolytopes}
states that the regular triangulations $\Triang$ of $P$ (with vertices in $Z$)
are in correspondence with (reduced) Gr\"obner bases of $I_P$. The heights
inducing the triangulation yield a term order on $\C[ x_{\p} : \p \in Z]$ and
we write $\underline{x^R} - x^S$ to emphasize that $x^R$ is the leading term.
We set the following result on record, which reflects the content of
Theorem~8.3 and Corollary~8.9 of \cite{Sturmfel-GrobnerPolytopes}.  For
details on this algebraic-geometric correspondence outlined above, we recommend
the very accessible Chapter 8 of Sturmfels'
book~\cite{Sturmfel-GrobnerPolytopes}. 

\begin{theorem}\label{thm:sturmfels}
    Let $P \subset \R^d$ be a lattice polytope and let $\Triang$ be a
    regular and unimodular triangulation. Then a reduced Gr\"obner basis of
    $I_P$ is given by the collection of monomials
    \[
        \underline{x^R} - x^S \, ,
    \]
    where $R \subset P \cap \Z^d$ is a minimal non-face, $S$ is a
    multisubset of $P \cap \Z^d$ such that $\sum R = \sum S$ and $\conv(S)$
    is a face of some simplex in $\Triang$.
    In particular, $\Triang$ is flag if and only if the leading terms are
    quadratic and square-free. 
\end{theorem}

Let $\Triang$ be a unimodular triangulation of $P$.  Given any lattice point
$\p \in 2 P$ there are unique $\p^{(1)}, \p^{(2)} \in Z$ such that $\p =
\p^{(1)} + \p^{(2)}$ and $\conv(\p^{(1)},\p^{(2)})$ is a face of a simplex
in $\Triang$. Note that $\p^{(1)}$ and $\p^{(2)}$ do not have to be distinct
points.  Let us call two points $\p,\q \in \R^d$ \Def{separable} if $p_i$
and $q_i$ have different signs for some $i=1,\dots,d$.  Together with
Theorem~\ref{thm:ExtendTriang}, this yields the following description of a
Gr\"obner basis for unconditional reflexive polytopes.

\begin{theorem}\label{thm:grobner_general}
    Let $P \subset \R^d$ be an anti-blocking polytope with a regular and
    unimodular triangulation and let $\underline{x^{R_i}} - x^{S_i}$ for
    $i=1,\dots,m$ be the associated Gr\"obner basis for $I_{P}$. A Gr\"obner
    basis associated to the toric ideal of $\U P$ is given by the binomials:
    \[
        \underline{x^{\sigma R_i}} - x^{\sigma S_i}
    \]
    for $i=1,\dots,m$ and $\sigma \in \PM^d$ and 
    \[
        \underline{x^{\p} \, x^{\q}} \ - \ x^{\sigma \e^{(1)}} \, x^{\sigma
        \e^{(2)}}
    \]
    for any separable $\p,\q \in \U
    P\cap\Z^d$ and $\sigma \in \PM^d$ such that $\sigma(\p+\q) = \e \in 2 P$.
\end{theorem}
\begin{proof}
    Theorem~\ref{thm:ExtendTriang} states that the regular and unimodular
    triangulation $\Triang$ of $P$ induces a regular and unimodular
    triangulation $\U \Triang$ of $\U P$.  It follows from
    Theorem~\ref{thm:ExtendTriang} (iii) that a minimal non-face $R$ of
    $\U\Triang$ is of the form $R = \sigma R'$, where $R'$ is a non-face of
    $\Triang$ and $\sigma \in \PM^d$, or it is of the form $R = \{\p,\q\}$ for
    separable $\p,\q \in Z := \U P \cap \Z^d$.

    In order to apply Theorem~\ref{thm:sturmfels}, we need to determine for
    every minimal non-face $R$ the multisubset $S$ of $Z$ such that $x^R - x^S
    \in I_{\U P}$.
    If $R = \sigma
    R_i$ for some minimal non-face $R_i$, then we can take $S = \sigma S_i$.
    If $R = \{ \p,\q \}$, then there is some $\sigma \in \PM^d$, such that
    $\sigma(\p + \q) = \e \in 2P$. It follows that $\p + \q \in \conv( 
    \sigma \e^{(1)}, \sigma \e^{(2)})$, which is a face of some simplex of
    $\U\Triang$. Hence we can take $S = \{ \sigma \e^{(1)}, \sigma \e^{(2)}
    \}$.
\end{proof}

A prominent class of perfect graphs $G$ for which regular, unimodular
triangulations of  $P_G$, as well as Gr\"obner bases for $I_{P_G}$, are well
understood are comparability graphs of finite posets.  Let $\Pi =
([d],\preceq)$ be a partially ordered set with comparability graph
$G_{\prec}$.  The stable set polytopes of comparability graphs were studied
by Stanley~\cite{StanleyOrder} under the name \Def{chain polytopes} and are
denoted by $\Chain(\Pi)$.  An \Def{antichain} in $\Pi$ is a collection of
pairwise incomparable elements.  The vertices of $P_{G_{\prec}}$ are
precisely the points $\1_A$, where $A$ is an antichain.  Let $\Anti(\Pi)$
denote the collection of antichains.  A pulling triangulation of
$P_{G_\prec}$ can be explicitly described (see Section 4.1
in~\cite{ChappellFriedlSanyal} for exposition and details). The
corresponding (reverse lexicographic) Gr\"obner basis was described by
Hibi~\cite{Hibi-LatticesAndStraighteningLaws}.
Following~\cite{ChappellFriedlSanyal}, we define 
\[
    A \sqcup A' \ \coloneqq \ \min(A \cup A') \quad \text{ and } \quad 
    A \sqcap A' \ \coloneqq  \ (A \cap A') \cup (\max(A \cup A') \setminus \min(A \cup A')) \, ,
\]
where min and max are taken with respect to the partial order $\preceq$.  We
call two antichains $A,A'$ \Def{incomparable} if $\max(A \cup A')$ is not a
subset of $A$ and not of $A'$. To ease notation, we identify variables
$x_A$ in $\C[ x_A : A \in \Anti(\Pi)]$ with symbols $[A]$.

\begin{theorem}[\cite{Hibi-LatticesAndStraighteningLaws,StanleyOrder}]\label{thm:grobner_chain}
    Let $\Pi$ be a poset and $\Chain(\Pi)$ its chain polytope.  A Gr\"obner
    basis for $I_{\Chain(\Pi)}$ is given by the binomials
    \[
        \underline{[B] \cdot [B']} - [B \sqcup B'] \cdot [B \sqcap B'] \, ,
    \]
    for all incomparable antichains $B,B' \in \Anti(\Pi)$. The corresponding
    triangulation of  $\Chain(\Pi)$ is regular, unimodular, and flag. 
\end{theorem}

We define the \Def{unconditional chain polytope} $\UC(\Pi)$ as the
unconditional reflexive polytope associated to $G_\prec$. These polytopes
were independently introduced by Ohsugi and Tsuchiya~\cite{OT18} under the
name \emph{enriched chain polytopes}.  The lattice points in $\UC(\Pi)$ are
uniquely given by 
\[
    \1_B - 2\cdot \1_A \, ,
\]    
where $A \subseteq B$ are antichains. We write $B - A$ for the pair $A
\subseteq B$ of antichains. In the following, we slightly abuse notation and
write $B - C$ instead of $B - (B \cap C)$ for anti-chains $B,C$. 
We get a description of the vertices of $\UC(\Pi)$ from
Proposition~\ref{prop:AB-inequalities}: Every vertex of $\UC(\Pi)$ is of the
form $\vb = \sigma \1_B$ for some inclusion-maximal antichain $B \subseteq
\Pi$.  Setting $A := \{ i \in B : v_i < 0 \}$, we deduce that the vertices of
$\UC(\Pi)$ are uniquely given by $B - A$ where $B$ is an inclusion-maximal
anti-chain.

We call $B-A$ and $B'-A'$ separable if the corresponding points are
separable. We also extend the two operations introduced above
\begin{align*}
    (B - A) \sqcup (B' - A') \ &:= \ (B \sqcup B') - (A \sqcup A')  \\
    (B - A) \sqcap (B' - A') \ &:= \ (B \sqcap B') - (A \sqcap A') 
\end{align*}

The following result was also obtained in~\cite{OT18}.

\begin{theorem}\label{thm:grobner_uchain}
    Let $\Pi = ([d],\preceq)$ be a finite poset and $I_{\UC(\Pi)}$ the toric
    ideal associated to the unconditional chain polytope $\UC(\Pi)$. Then a
    reduced Gr\"obner basis is given by the binomials 
    \[
        \underline{[B-E] \cdot [B'-E]} \ - \ 
        [(B \sqcup B')-E] \cdot [(B \sqcap B')-E] \, ,
    \]
    for all incomparable $B,B' \in \Anti(\Pi)$ and $E \subseteq B \cup B'$ as
    well as 
    \[
        \underline{[B - A] \cdot [B' - A']} \ - \ 
        [ (C - D) \sqcup (C' - D') ] \cdot
        [ (C - D) \sqcap (C' - D') ] \, 
    \]
    where $B-A, B'-A'$ are separable and
    \begin{align*}
    C  &:= (B \setminus B') \cup (A \cap A') &
    D  &:= (A \setminus B') \cup (A \cap A') \\
    C' &:= (B' \setminus B) \cup (A \cap A') &
    D' &:= (A' \setminus B) \cup (A \cap A')  \, .
    \end{align*}
\end{theorem}
\begin{proof}
    We apply Theorems~\ref{thm:grobner_general} using the Gr\"obner basis of
    the toric ideal of $\Chain(\Pi)$ provided by
    Theorem~\ref{thm:grobner_chain}. This provides the first set of binomials
    and we only have to argue the binomials with leading terms coming from
    separable pairs $B-A$ and $B'-A'$.

    The Gr\"obner basis given in Theorem~\ref{thm:grobner_chain} tells us how
    to find the decomposition $\e^{(1)} + \e^{(2)}$ of a point $\e \in 2
    \Chain(\Pi)$. Indeed, the minimal non-faces of the triangulation are given
    by the leading terms and correspond to incomparable pairs of anti-chains
    $E,E' \in \Anti(\Pi)$. The point $\1_E + \1_{E'}$ can then be written as
    $\1_{E \sqcup E'} + \1_{E \sqcap E'}$. The anti-chains $E \sqcup E', E
    \sqcap E'$ are comparable and thus $\conv(\1_{E \sqcup E'} , \1_{E \sqcap
    E'})$ is a face a simplex in the triangulation.

    By our discussion above, every point in $\U \Chain(\Pi)$ is of the form
    $\1_B - 2 \cdot \1_A$ for $B \in \Anti(\Pi)$ and $A \subseteq B$. Now, for $B-A$
    and $B'-A'$ separable, we have
    \[
        \p \ = \ (\1_B - 2 \cdot \1_A) + (\1_{B'} - 2\cdot \1_{A'})  \ = \ 
        (\1_{C} - 2 \cdot \1_{D}) + 
        (\1_{C'} - 2 \cdot \1_{D'}) \, .
    \]
    Note that $C - D$ and $C' - D'$ are not separable and the pairs $(C - D)
    \sqcup (C' - D')$ and $(C - D) \sqcap (C' - D')$ give a suitable
    representation of $\p$ with respect to the triangulation of $\U
    \Chain(\Pi)$.
\end{proof}

\section{Concluding remarks} 
\label{sec:conclude}

\subsection{A Blaschke--Santal\'o inequality}
The Blaschke--Santal\'o inequality~\cite[Sect.~10.7]{Schneider} implies that
for a centrally-symmetric convex body $K \subset \R^d$
\[
    \vol(K) \cdot \vol(K^\ast) \ \le \ \vol(B_d)^d,
\]
where $B_d$ is the Euclidean unit ball. Equality is attained precisely when
$K$ is an ellipsoid.  Based on computations for up to $9$ vertices, we
conjecture the following.

\begin{conjecture}
    For every $n \ge 1$, there is a unique perfect graph $G$ on $n$ vertices
    such that $\vol(\U P_G) \cdot \vol( \U P_{\overline{G}})$ is maximal.
\end{conjecture}

For $n=3,4,5$, the unique maximizer is the path on $n$ vertices. For $n=6$ the
unique maximizer is a $6$-cycle and for $n=7$ the maximizer is obtained by
adding a dangling edge to the $6$-cycle. For $n=8$ the graph in question is
$6$-cycle with an additional $4$-path connecting two antipodal vertices but
for $n=9$ the graph is much more complicated.

\subsection{Birkhoff polytopes of other types}
It is only natural to look at Birkhoff-type polytopes of other finite
irreducible Coxeter groups.  Since the type-$B$ and the type-$C$ Coxeter
groups are equal, we get the same polytope.  Recall that the type-$D$ Coxeter
group $D_n$ is the subgroup of $B_n$ with permutations with an even number of
negative entries. We can construct the \Def{type-$D$ Birkhoff polytope},
$\mathcal{BD}(n)$, to be the convex hull of signed permutation matrices with
an even number of negative entries.  As one may suspect from this
construction, the omission of all lattice points in various orthants which
occurs in $\mathcal{BD}(n)$ ensures that it cannot be an OLP polytope and is
thus not subject to any of our general theorems.  When $n=2$ and $n=3$,
$\mathcal{BD}(n)$ is a reflexive polytope, but $\mathcal{BD}(3)$ does not have
the IDP.  Moreover, $\mathcal{BD}(4)$ fails to be reflexive.

Additionally, one could consider Birkhoff constructions for Coxeter groups of exceptional type, in particular $E_6$, $E_7$ and $E_8$ (see, e.g., \cite{CombinatoricsCoxeter}).
While we did not consider these polytopes in our investigation, we do raise the following question:

\begin{question}
Do the Birkhoff polytope constructions for $E_6$, $E_7$, and $E_8$ have the IDP?
Are these polytopes reflexive? 
Do they have other interesting properties?
\end{question}

\subsection{Future directions}\label{sec:future}
In addition to considering Birkhoff polytopes of other types and connections
to Gale duality as discussed above, there are several  immediate avenues for
further research.  Coxeter groups are of great interest in the broader
community of algebraic and geometric combinatorics (see, e.g.,
\cite{CombinatoricsCoxeter}) and it would be interesting if Coxeter-theoretic
insights can be gained from the geometry of $\BB(n)$.

An additional future direction is to consider applications of the
orthant-lattice property, particularly those of Theorem~\ref{thm:ExtendTriang}
and Remark \ref{rmk:generalIDP}.  One potentially fruitful avenue is an
application to reflexive \emph{smooth} polytopes.  Recall that a lattice
polytope $P\subset \R^d$ is \Def{simple} if every vertex of $P$ is contained
in exactly $d$ edges (see, e.g., \cite{Ziegler}).  A simple polytope $P$ is
called \Def{smooth} if the primitive edge direction generate $\Z^d$ at every
vertex of $P$.  Smooth polytopes are particularly of interest due to a
conjecture commonly attributed to Oda~\cite{OdasConjecture}:
\begin{conjecture}[Oda]
    If $P$ is a smooth polytope, then $P$ has the IDP.
\end{conjecture}

This conjecture is not only of interest in the context of Ehrhart theory, but
also in toric geometry.  One potential strategy is to consider similar
constructions to OLP polytopes for smooth reflexive polytopes to make progress
towards this problem.  As a first step, we pose the following question:

\begin{question}
    Are all smooth reflexive polytopes OLP polytopes? 
\end{question}

Furthermore, regarding reflexive OLP polytopes one can ask the question:

\begin{question}
    Given a reflexive OLP polytope $P$, under what conditions can we
    guarantee that $P^\ast$ is a reflexive OLP polytope?
\end{question}

By~\eqref{eqn:complement}, this has a positive answer when
$P$ is an unconditional reflexive polytope.  However, there are
multiple examples of failure in general even in dimension $2$ (see
Figure~\ref{16reflexivePolygons-Picture}).

\bibliographystyle{amsalpha}
\bibliography{UnconditionalReflexive}
\end{document}